\documentclass[a4paper,11pt]{article}
\usepackage{amsmath}
\usepackage{dsfont}
\usepackage{amssymb}
 \usepackage[usenames,dvipsnames]{pstricks}
 \usepackage{epsfig}
 \usepackage{pst-grad} 
 \usepackage{pst-plot} 
\numberwithin{equation}{section}
 \title{Estimates on the speedup and slowdown for a diffusion in a drifted brownian potential.}
 \author{Gabriel Faraud\footnote{Université Paris 13
CNRS, UMR 7539 LAGA
99, avenue Jean-Baptiste Clément
F-93 430 Villetaneuse,
France research partially supported by the ANR project MEMEMO.}.}
  \newtheorem{theorem}{Theorem}[section]

\newtheorem{lemma}[theorem]{Lemma}
\newtheorem{statement}{Statement}[section]
\DeclareMathOperator{\egall}{\overset{law}{=}}
 
\begin{document}
 \begin{titlepage}\maketitle
 \begin{center}
 \date
 {\tt faraud@math.univ-paris13.fr}
 
 116 rue de Lagny 
 
 93100 Montreuil
 
 France
 
 (+33)149404085
 \end{center}
\begin{abstract}
We study a model of diffusion in a brownian potential. This model was firstly introduced by T. Brox (1986) as a continuous time analogue of random walk in random environment. We estimate the deviations of this process above or under its typical behavior. Our results rely on different tools such as a representation introduced by Y. Hu, Z. Shi and M. Yor, Kotani's lemma, introduced at first by K. Kawazu and H. Tanaka (1997), and a decomposition of hitting times developed in a recent article by A. Fribergh, N. Gantert and S. Popov (2008) \nocite{Gantert:2008}. Our results are in agreement with their results in the discrete case.
\end{abstract}
\end{titlepage}
\section{Introduction.} 
Process in random media have been introduced in order to study physical or biological mechanisms such as the replication of DNA. The first model, in discrete time, goes back to A. Chernov \cite{chernov1967rmc} and D. Temkin \cite{temkin1972odr}. It is now well understood : see, for example \cite{solomon1975rwr}, \cite{sinai1983lbo}, or \cite{kesten1975llr}.
A continuous time version of this process has been introduced by S. Schumacher \cite{schumacher1985drc}, and studied by T. Brox \cite{brox:1986}. It can be described as follows.

Let $(W(x))_{x\in \mathbb{R}}$ be a one-dimensional brownian motion
 defined on $\mathbb{R}$ starting from $0$, and, for $\kappa\in
  \mathbb{R}$, $$W_\kappa(x):=W(x)-\frac{\kappa}{2}x.$$

 Let $(\beta(t))_{t\geq 0}$
 be another one-dimensional brownian motion, independent of $W$. 
We call {\it diffusion process with potential $W_{\kappa}$} a solution to the (formal) equation 
\begin{equation}\label{eds}dX_t=d\beta_t -\frac{1}{2} W_{\kappa}'(X_t)dt.\end{equation}
$W_{\kappa}'$ has clearly no rigorous meaning, but a mathematical definition of (\ref{eds}) can be given in terms of the
infinitesimal generator. For a given realization of $W_\kappa$, $X_t$ is
a real-valued diffusion started at $0$ with generator 
$$\frac{1}{2} e^{W_\kappa(x)}\frac{d}{d x}\left(e^{-W_\kappa(x)} \frac{d}{d x}\right).$$

This diffusion can also be defined by a time-change representation : 
$$X_t=A_\kappa^{-1}\left(B(T_\kappa^{-1}(t))\right),$$ where
$$A_\kappa(x)=\int_0^x e^{W\kappa(y)}d y,$$
$$T_\kappa(t)=\int_0^t e^{-2W_\kappa(A_\kappa^{-1}(B(s)))}d s,$$
and $B$ is a standard Brownian motion.
$A_\kappa$ is the scale function of this process, and its speed measure is $2e^{-W_{\kappa}(x)}d x$.

\begin{figure}
\scalebox{1} 
{
\begin{pspicture}(-2,-2.7239997)(8.204,2.76)
\rput(3.128,0.46){\psaxes[linewidth=0.02,labels=none,ticks=none,ticksize=0.10583333cm]{->}(0,0)(-3,-3)(5,2)}
\pscustom[linewidth=0.0080]
{
\newpath
\moveto(0.0,2.3122287)
\lineto(0.026282044,2.3317332)
\curveto(0.039423063,2.3414857)(0.07884613,2.3707428)(0.10512818,2.3902476)
\curveto(0.13141021,2.4097524)(0.16645302,2.4146285)(0.1752137,2.4)
\curveto(0.18397439,2.3853714)(0.20149574,2.3561144)(0.21025643,2.3414857)
\curveto(0.2190171,2.326857)(0.23653847,2.2829714)(0.24529915,2.2537143)
\curveto(0.25405982,2.224457)(0.27596152,2.1952)(0.28910255,2.1952)
\curveto(0.3022436,2.1952)(0.33290604,2.2147048)(0.3504274,2.2342095)
\curveto(0.36794877,2.2537143)(0.40299147,2.2585905)(0.42051286,2.243962)
\curveto(0.4380342,2.2293334)(0.47745728,2.1903238)(0.49935898,2.165943)
\curveto(0.5212607,2.141562)(0.5694444,2.1269333)(0.5957265,2.1366856)
\curveto(0.6220085,2.1464381)(0.665812,2.0976763)(0.68333334,2.039162)
\curveto(0.7008547,1.9806476)(0.74027777,1.9123809)(0.7621795,1.9026285)
\curveto(0.7840812,1.8928761)(0.82350427,1.8782476)(0.8410256,1.8733715)
\curveto(0.858547,1.8684952)(0.88920945,1.8099809)(0.9023505,1.7563429)
\curveto(0.91549146,1.7027048)(0.9373932,1.6100571)(0.9461539,1.5710477)
\curveto(0.95491457,1.5320381)(0.9943376,1.4881525)(1.025,1.4832762)
\curveto(1.0556624,1.4784)(1.1169871,1.5320381)(1.1476495,1.5905523)
\curveto(1.178312,1.6490666)(1.2264956,1.7319618)(1.244017,1.7563429)
\curveto(1.2615384,1.7807238)(1.3097222,1.8489904)(1.3403845,1.8928761)
\curveto(1.371047,1.9367619)(1.4104701,1.9611428)(1.4192308,1.941638)
\curveto(1.4279915,1.9221333)(1.4542736,1.8587428)(1.4717948,1.8148571)
\curveto(1.4893162,1.7709714)(1.511218,1.7027048)(1.5155983,1.6783239)
\curveto(1.5199786,1.6539428)(1.5331197,1.6003048)(1.5418804,1.5710477)
\curveto(1.5506411,1.5417905)(1.576923,1.4881523)(1.5944444,1.4637713)
\curveto(1.6119658,1.4393904)(1.6426282,1.4247618)(1.6557692,1.4345142)
\curveto(1.6689103,1.4442667)(1.7170941,1.483276)(1.7521368,1.5125333)
\curveto(1.7871796,1.5417905)(1.8353633,1.551543)(1.8485043,1.5320381)
\curveto(1.8616453,1.5125333)(1.9010684,1.4393905)(1.9273505,1.3857524)
\curveto(1.9536325,1.3321142)(2.0061965,1.2882286)(2.0324786,1.2979809)
\curveto(2.0587606,1.3077333)(2.1069443,1.3077333)(2.1288462,1.2979809)
\curveto(2.1507478,1.2882286)(2.190171,1.268724)(2.2076921,1.2589716)
\curveto(2.2252135,1.2492191)(2.2646368,1.2345904)(2.2865384,1.2297144)
\curveto(2.30844,1.2248381)(2.3391025,1.2004571)(2.3478632,1.1809524)
\curveto(2.356624,1.1614478)(2.369765,1.1126858)(2.3741453,1.0834286)
\curveto(2.3785255,1.0541714)(2.3916667,1.0102856)(2.4004273,0.995657)
\curveto(2.409188,0.98102844)(2.4398506,0.9615238)(2.461752,0.95664763)
\curveto(2.4836535,0.95177156)(2.5274572,0.93714297)(2.549359,0.92739046)
\curveto(2.571261,0.91763806)(2.6019232,0.8835048)(2.610684,0.8591238)
\curveto(2.6194446,0.8347429)(2.636966,0.785981)(2.6457267,0.7616)
\curveto(2.6544874,0.7372191)(2.6807694,0.67870486)(2.6982908,0.64457154)
\curveto(2.715812,0.6104382)(2.742094,0.53241915)(2.7508547,0.48853332)
\curveto(2.7596154,0.44464752)(2.768376,0.37638092)(2.768376,0.35200012)
\curveto(2.768376,0.32761902)(2.7858975,0.28860962)(2.8034189,0.273981)
\curveto(2.8209403,0.25935242)(2.8603632,0.2544763)(2.882265,0.26422852)
\curveto(2.904167,0.2739807)(2.9304485,0.30811402)(2.934829,0.33249512)
\curveto(2.9392097,0.35687622)(2.9567308,0.3910095)(2.9698718,0.40076202)
\curveto(2.9830127,0.41051438)(3.0136752,0.43001908)(3.0311966,0.4397714)
\curveto(3.048718,0.44952378)(3.096902,0.4592761)(3.1275642,0.4592761)
\curveto(3.1582263,0.4592761)(3.2107904,0.43489516)(3.2326922,0.41051424)
\curveto(3.254594,0.38613328)(3.3159187,0.31299043)(3.355342,0.26422852)
\curveto(3.394765,0.21546662)(3.4736109,0.098438114)(3.513034,0.030171508)
\curveto(3.5524573,-0.038095094)(3.6137824,-0.1404953)(3.635684,-0.1746286)
\curveto(3.6575854,-0.2087619)(3.7101495,-0.2721521)(3.740812,-0.3014093)
\curveto(3.7714746,-0.3306665)(3.8240387,-0.394057)(3.8459404,-0.42819032)
\curveto(3.8678417,-0.4623236)(3.8985043,-0.5159619)(3.907265,-0.5354669)
\curveto(3.9160256,-0.55497164)(3.9423077,-0.6086096)(3.959829,-0.64274293)
\curveto(3.9773505,-0.67687625)(4.021154,-0.7353906)(4.0474358,-0.7597714)
\curveto(4.073718,-0.7841522)(4.121902,-0.83779055)(4.1438036,-0.8670477)
\curveto(4.1657047,-0.8963049)(4.2226496,-0.94019043)(4.2576923,-0.954819)
\curveto(4.292735,-0.9694476)(4.384722,-1.0133331)(4.4416666,-1.0425904)
\curveto(4.498611,-1.0718476)(4.5818377,-1.0962286)(4.6081195,-1.0913526)
\curveto(4.634402,-1.0864764)(4.686966,-1.0523431)(4.713248,-1.023086)
\curveto(4.73953,-0.9938287)(4.774573,-0.93531436)(4.7833333,-0.9060571)
\curveto(4.7920933,-0.87679994)(4.8139954,-0.7939047)(4.827136,-0.74026674)
\curveto(4.840277,-0.6866284)(4.86656,-0.6086096)(4.879701,-0.5842285)
\curveto(4.892842,-0.5598474)(4.9278846,-0.50620973)(4.949786,-0.47695252)
\curveto(4.9716883,-0.4476953)(5.0461535,-0.394057)(5.098718,-0.3696762)
\curveto(5.1512823,-0.3452954)(5.238889,-0.3209143)(5.273932,-0.3209143)
\curveto(5.3089747,-0.3209143)(5.365919,-0.3696762)(5.387821,-0.4184381)
\curveto(5.409723,-0.4672)(5.4360046,-0.55497134)(5.4403844,-0.593981)
\curveto(5.444765,-0.6329907)(5.4666667,-0.69150484)(5.484188,-0.7110095)
\curveto(5.5017095,-0.7305142)(5.541133,-0.7500189)(5.5630345,-0.7500189)
\curveto(5.5849366,-0.7500189)(5.6199794,-0.7353903)(5.63312,-0.7207617)
\curveto(5.646261,-0.7061334)(5.6813035,-0.71588564)(5.7032056,-0.74026674)
\curveto(5.725107,-0.76464784)(5.76453,-0.8231619)(5.7820516,-0.8572952)
\curveto(5.799573,-0.89142853)(5.8871794,-0.94019043)(5.957265,-0.954819)
\curveto(6.0273504,-0.9694476)(6.163141,-0.9548193)(6.2288465,-0.92556214)
\curveto(6.2945514,-0.8963049)(6.395299,-0.8963049)(6.4303417,-0.92556214)
\curveto(6.4653845,-0.9548193)(6.513568,-1.0133331)(6.526709,-1.0425904)
\curveto(6.5398498,-1.0718476)(6.6055555,-1.0962286)(6.6581197,-1.0913526)
\curveto(6.710684,-1.0864764)(6.78953,-1.1059809)(6.8158116,-1.1303619)
\curveto(6.8420935,-1.1547431)(6.885897,-1.2278857)(6.9034185,-1.2766477)
\curveto(6.92094,-1.3254095)(6.9910254,-1.4814477)(7.0435896,-1.5887238)
\curveto(7.0961537,-1.6959997)(7.188141,-1.8715429)(7.227564,-1.9398096)
\curveto(7.2669873,-2.0080762)(7.3589745,-2.1446095)(7.4115386,-2.212876)
\curveto(7.4641027,-2.2811425)(7.547329,-2.383543)(7.577992,-2.4176764)
\curveto(7.608654,-2.4518096)(7.709402,-2.529829)(7.779487,-2.5737145)
\curveto(7.8495727,-2.6176002)(7.9766026,-2.6712377)(8.033547,-2.68099)
\curveto(8.090492,-2.6907423)(8.160578,-2.7053711)(8.2,-2.72)
}

\pscustom[linewidth=0.000]
{
\newpath
\moveto(3.82,-0.32)
 \curveto(3.7714746,-0.3306665)(3.8240387,-0.394057)(3.8459404,-0.42819032)
\curveto(3.8678417,-0.4623236)(3.8985043,-0.5159619)(3.907265,-0.5354669)
\curveto(3.9160256,-0.55497164)(3.9423077,-0.6086096)(3.959829,-0.64274293)
\curveto(3.9773505,-0.67687625)(4.021154,-0.7353906)(4.0474358,-0.7597714)
\curveto(4.073718,-0.7841522)(4.121902,-0.83779055)(4.1438036,-0.8670477)
\curveto(4.1657047,-0.8963049)(4.2226496,-0.94019043)(4.2576923,-0.954819)
\curveto(4.292735,-0.9694476)(4.384722,-1.0133331)(4.4416666,-1.0425904)
\curveto(4.498611,-1.0718476)(4.5818377,-1.0962286)(4.6081195,-1.0913526)
\curveto(4.634402,-1.0864764)(4.686966,-1.0523431)(4.713248,-1.023086)
\curveto(4.73953,-0.9938287)(4.774573,-0.93531436)(4.7833333,-0.9060571)
\curveto(4.7920933,-0.87679994)(4.8139954,-0.7939047)(4.827136,-0.74026674)
\curveto(4.840277,-0.6866284)(4.86656,-0.6086096)(4.879701,-0.5842285)
\curveto(4.892842,-0.5598474)(4.9278846,-0.50620973)(4.949786,-0.47695252)
\curveto(4.9716883,-0.4476953)(5.0461535,-0.394057)(5.098718,-0.3696762)
\curveto(5.1512823,-0.3452954)(5.238889,-0.3209143)(5.273932,-0.3209143)
\gsave
 \psline[liftpen=1](5.22,-0.32)(3.82,-0.32)
 \fill[fillstyle=solid,fillcolor=gray]
 \grestore 
}

\psline[linewidth=0.0080cm,linestyle=dashed,dash=0.16cm 0.16cm](5.22,-0.32)(3.82,-0.32)

\rput(5.31,0.1){\footnotesize Valley}
\psline[linewidth=0.0080cm,arrowsize=0.05291667cm 2.0,arrowlength=1.4,arrowinset=0.4]{->}(4.74,0.08)(4.5,-0.58)

\rput(3.81,2.56){\footnotesize $W_{\kappa}(x)$}

\rput(8.13,0.105){x}
\psline[linewidth=0.022cm,arrowsize=0.05291667cm 2.0,arrowlength=1.4,arrowinset=0.4]{->}(4.88,1.32)(6.34,1.32)

\rput(5.59,1.625){\footnotesize drift} 
\end{pspicture} 
}
\label{potentieldessin}
\caption{ A ``valley''.}
\end{figure}
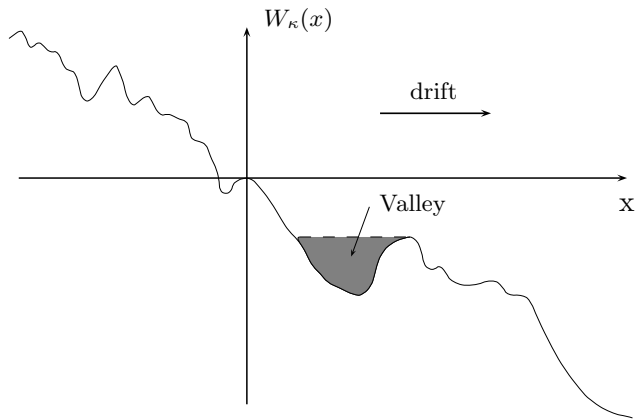

Intuitively, for a given environment $W_{\kappa}$, the diffusion $X_t$ will tend to go to places where $W_{\kappa}$ is low, and to spend a lot of time in the ``valleys'' of $W_{\kappa}$. If the environment is drifted $(\kappa>0)$, the process will be transient to the right, but it will be slowed by those valleys (see figure \ref{potentieldessin}). This will be explained more precisely in section \ref{quenchedslowdown}.

For general background on diffusion processes and time-change representation we refer to \cite{Rogers:1986,Revuz:1994,ito1974dpa}.

We will call $\mathcal{P}$ the probability associated to $W$, $P_W$
the quenched probability associated to the diffusion, and
$\mathbb{P}:=\mathcal{P}\otimes P_W$ the annealed probability.

T. Brox gave a result concerning the long time behavior of the diffusion in the case $\kappa=0$. Namely, under the probability $\mathbb{P},$
$$\frac{X_t}{(\log t)^2}\rightarrow U,$$ where $U$ follows an explicit distribution.

 The case $\kappa>0$ was studied both by K. Kawazu and H. Tanaka (\cite{Kawazu:1997})
and Y. Hu, Z. Shi, M. Yor (\cite{Hu:1999}) and exhibits a
``Kesten-Kozlov-Spitzer'' behavior: when $\kappa>1$, the diffusion has a positive speed; when $\kappa=1$, under $\mathbb{P},$
$$\frac{X_t \log t}{t}\rightarrow 4$$ in probability, while, when $0<\kappa<1,$ 
$$\frac{X_t}{t^\kappa}\rightarrow V$$ in distribution, where $V$ follows the inverse of a completely asymmetric stable law.

We are interested in the deviations between $X_t$ and its asymptotic behavior, in the case $0<\kappa<1$.

This questions have already been studied in the other cases, we refer to \cite{Hu:2004} for estimates in the case $\kappa=0$, and to \cite{Talet:2007} for large deviation estimates in the case $\kappa>1$.
 \\

Our study will split into four different problems, indeed the quenched and annealed settings present different behavior, and for each of them we have to consider deviations above the asymptotic behavior (or speedup) and deviations under the asymptotic behavior (or slowdown). 

We start with the annealed results.
 For $u$ and $v$ two
functions of $t$, we note $u\gg v$ if $u/v\rightarrow_{t\rightarrow \infty} \infty.$
\begin{theorem}[Annealed speedup/slowdown]\label{as}
Suppose $0<\kappa <1$, and $u\rightarrow \infty$ is a function of $t$ such that for some $\varepsilon>0$, $u\ll t^{1-\kappa-\varepsilon}$,
then there exist two positive constants $C_1$ and $C_2$ such that
\begin{equation}\label{ase1}\lim_{t\rightarrow \infty}
\frac{-\log \mathbb{P}\left(X_t>t^\kappa u \right)}{u^{\frac{1}{1-\kappa}}}=C_1,\end{equation}
and if $\log{u}\ll t^\kappa$, 
\begin{equation} \label{ase2}\lim_{t\rightarrow \infty} u\mathbb{P}\left(X_t<\frac{t^\kappa}{u} \right)=C_2.\end{equation} 
Furthermore, the results remain true if we replace $X_t$ by $\sup_{s<t} X_s$.
\end{theorem}

This is in fact a easy consequence of the study of the hitting time of a certain level by the diffusion.
We set $H(v) = \inf\{t>0: X_t=v\}.$
We have the following estimates.
\begin{theorem}\label{pas} Suppose $0<\kappa <1$ and $\varepsilon>0$.
For $u\rightarrow \infty$ $v\rightarrow \infty$ two functions of $t$ such that for some $\varepsilon>0$, $u\ll v^{1-\kappa-\varepsilon}$, there exist two positive constants $C_1$ and $C_2$ such that 
\begin{equation}\label{htasu}\lim_{t\rightarrow \infty}\frac{-\log\mathbb{P}\left[H(v)<\left(\frac{v}{u}\right)^{1/\kappa}\right]}{u^{\frac{1}{1-\kappa}}}=C_1,\end{equation}
and if $\log{u}\ll v$, 
\begin{equation}\label{htasd} \lim_{t\rightarrow \infty} u\mathbb{P}\left[H(v)>(v u)^{1/\kappa}
 \right]=C_2. \end{equation} 
\end{theorem}
The proof of this result involves a representation of $H(v)$ introduced in \cite{Hu:2004}. 
\\
 
We now turn to the quenched setting. We have the following estimates for the speedup
\begin{theorem}[Quenched speedup]\label{qs}
Suppose $0<\kappa <1$, and $u\rightarrow \infty$ is a function of $t$ such that for some $\varepsilon>0$, $u\ll t^{1-\kappa-\varepsilon}$,
then there exists a positive constants $C_3$  such that

$$\lim_{t\rightarrow\infty}\frac{-\log P_{W}\left(X_t>t^\kappa u \right)}{u^{\frac{1}{1-\kappa}}}= C_3,\;P-a.s..$$
Furthermore the result remains true if we replace $X_t$ by $\sup_{s\leq t} X_s$.
\end{theorem}

As before, the proof of this will reduce to estimates on the hitting times.
\begin{theorem}	\label{pqs}	
For $u\rightarrow \infty$ $v\rightarrow \infty$ two functions of $t$ such that for some $\varepsilon>0$, $u\ll v^{1-\kappa-\varepsilon},$ then
\begin{equation} \lim_{t\rightarrow\infty}\frac{-\log P_W\left[H(v)<\left(\frac{v}{u}\right)^{1/\kappa}\right]}{u^{\frac{1}{1-\kappa}}}=C_3,\; P-a.s..\end{equation}
\end{theorem} 

For the slowdown, our result is less precise.
\begin{theorem}[Quenched slowdown]\label{pqsl}
Suppose $\kappa>0$. Let $\nu\in(0,1\wedge\kappa),$ then
\begin{equation}\label{pqsl1}
\lim_{t\rightarrow \infty}\frac{\log (-\log P_{W}[H(t^\nu)>t])}{\log t} =\left(1-\frac{\nu}{\kappa}\right)\wedge \frac{\kappa}{\kappa+1},\; P-a.s.,
\end{equation}
\begin{equation}\label{pqsl2}
\lim_{t\rightarrow \infty}\frac{\log (-\log P_{W}[X_{t}<t^\nu])}{\log t} =\left(1-\frac{\nu}{\kappa}\right)\wedge \frac{\kappa}{\kappa+1},\; P-a.s..
\end{equation}
\end{theorem}

Corresponding results for Random Walk in Random Environment have been developed in a recent article from A. Fribergh, N. Gantert and S. Popov \cite{Gantert:2008}. Our proof of the last result is quite inspired from theirs. 

The article will be organized as follows :
\begin{itemize}
\item In Section 2 we show Theorem \ref{as} and \ref{pas},
\item In Section 3 we show Theorem \ref{pqsl},
\item In Section 4 we show Theorem \ref{qs} and \ref{pqs}.
\end{itemize}
\section{The annealed estimate.}
For any nondecreasing function $u(t)$, we will denote by $u^{-1}
(t):=\inf\{v:u(v)>t\}$ the inverse function of $u$. We start with some preliminary statements.

\subsection{Preliminary statements.}

We first recall the Ray-Knight Theorems, they can be found in chapter XI of \cite{Revuz:1994}. Let $L_t^x$ be the local time at $x$ before $t$ of a brownian motion $\gamma_t$, and
  $\tau_t:={\left(L_.^0\right)}^{-1}(t)$ the inverse function of $L_t^0$. Let $\sigma(x)$ be the first
hitting time of $x$ by $\gamma_t$. 
\begin{statement}[First Ray-Knight Theorem]\label{rn1}
The process $\{L_{\sigma(a)}^{a-t}\}_{t\geq0}$ is a squared Bessel process, started at $0$, of dimension $2$ for $0\leq t\leq a$ and of dimension $0$ for $t\geq a$.
\end{statement}
\begin{statement}[Second Ray-Knight Theorem]\label{rn}Let $u\in\mathbb{R^+}$, 
The process $\{L_{\tau(u)}^t\}_{t\geq 0}$ is a squared Bessel process
of dimension $0$, starting from $u$.
\end{statement}

We have a useful representation of $H(v)$, due to Y. Hu and Z. Shi
(2004). Let $$\theta_1(v)=\int_0^{H(v)} \textbf{1}_{\{X_s\geq 0\}} d s,$$
and $$\theta_2(v)=\int_0^{H(v)} \textbf{1}_{\{X_s< 0\}} d s,$$
such that $H(v)=\theta_1(v)+\theta_2(v)$.
\begin{statement}\label{rep}
Let $\kappa\geq 0$ and $v>0$. Under $\mathbb{P}$, we have
$$\left(\theta_1(v),\theta_2(v)\right)\egall \left(4\int_0^v \left(
e^{\Xi_{\kappa} (s)}-1\right)d s, 16 \Upsilon
_{2-2\kappa}\left(e^{\Xi_{\kappa} (v) /2}\leadsto 1\right)\right).$$
Where $\Upsilon_{2-2\kappa} (x \leadsto y)$ denotes the first hitting
time of $y$ by a Bessel process of dimension $(2-2\kappa)$ starting from $x$, independent of
the diffusion $\Xi_{\kappa}$, which is the unique nonnegative solution
of
\begin{equation}\label{theta}
\Xi_{\kappa}(t) = \int_0^t \sqrt{1-e^{-\Xi_{\kappa}(s)}}d\beta'_s +
\int_0^t\left(-\frac{\kappa}{2}+\frac{1+\kappa}{2}e^{-\Xi_{\kappa}(s)}\right)ds,\;
t\geq0. 
\end{equation}
$\beta'$ being a standard brownian motion.
\end{statement}
We shall use the following lemma from \cite{Talet:2007}(Lemma 3.1).
\begin{statement}\label{mt}
Let $\{R_t\}_{t\geq 0}$ denote a squared Bessel process of dimension
$0$ started at $1$. For all $v,\delta>0$, we have  
$$P\left( \sup_{0\leq t\leq v} |R_t-1|>\delta \right) \leq
4\frac{\sqrt{(1+\delta)v}}{\delta}\exp{\left(-\frac{ \delta ^2}{8(1+\delta)v}\right)}.$$
\end{statement}

We now turn to the proof of Theorem \ref{pas}.
\subsection{Proof of Theorem \ref{pas}.}
Our proof will be separated in two parts : in the first part we will deal with the positive part of $H(v)$, $\theta_1$, then we will focus on $\theta_2$.
\subsubsection{The positive part.}
\label{tpp}
In view of statement \ref{rep}, we set
$$Z_t:= e^{\Xi_{\kappa} (t)}-1,$$ then $Z_t$ is the unique nonnegative
solution of 
$$d Z_t=\sqrt{Z_t(1+Z_t)} d\beta_t +\left(\frac{1-\kappa}{2} Z_t
  +\frac{1}{2} \right) d t,$$
and $$\theta_1(v)=4\int_0^v Z_t d t.$$ We call \begin{equation}\label{deff}f(z)=\int_1^z
  \frac{(1+s)^\kappa}{s}d s\end{equation} the scale function of $Z_t$.

We have $$f(Z_t)=\int_0^t \frac{(1+Z_s)^{\kappa
  +\frac{1}{2}}}{\sqrt{Z_s}} d\beta_s.$$
By the Dubbins-Schwartz representation (see chapter V, Theorem (1.6) of \cite{Revuz:1994}), there exists a standard
  Brownian motion $\gamma(t)$ such that
  \begin{equation} \label{eq}f(Z_t)=\gamma\left(\int_0^t \frac{(1+Z_s)^{2\kappa
  +1}}{Z_s}ds\right):=\gamma(\rho(t)).\end{equation}
We introduce \begin{multline}\label{defalpha}\alpha_t=\rho(t)^{-1}=\int_0^t
  \frac{Z_{\alpha_s}}{(1+Z_{\alpha_s})^{1+2\kappa}}d s \\= \int_0^t
  \frac{f^{-1}(\gamma_s)}{\left[1+ f^{-1} (\gamma_s)\right]^{1+
  2\kappa}}d s:=\int_0^t h(\gamma_s)d s.\end{multline}
We obtain easily the following equivalents
\begin{align*}
&f(z)\sim_{z\rightarrow \infty} z^\kappa /\kappa, \\
&f(z)\sim_{z\rightarrow 0} \log z, \\
&f^{-1}(z)\sim_{z\rightarrow \infty} (\kappa z)^{1/\kappa},\\
&f^{-1}(z)\sim_{z\rightarrow -\infty} e^z, \\
&h(z)\sim_{z\rightarrow \infty} (\kappa z)^{-2}, \\
&h(z)\sim_{z\rightarrow -\infty} e^z.
\end{align*}

We continue with a lemma, whose proof is postponed. Let $\tau_t$ be the inverse local time of $\gamma$.
\begin{lemma}\label{lemma} Let $\epsilon>0$, $c_h:=\int_0^\infty
  h(x)d x$. Let $w(t)\rightarrow \infty$, such that $w(t)/t\rightarrow0$. Then for $t$ large enough,
$$ \mathbb{P}\left(\rho(t)>\tau_{\frac{t}{(1-3\epsilon)c_h}}\right)\leq
  \exp\left(-{w}\right),$$
  and
$$\mathbb{P}\left(\rho(t)<\tau_{\frac{t}{(1+3\epsilon)c_h}}\right)\leq
  \exp\left(-{w}\right).$$
\end{lemma}

Let $\tilde{v}\ll v$, in view of (\ref{eq}),
\begin{equation}\label{defg}\theta_1(v)=4\int_0^v
f^{-1}(\gamma_{\rho(s)})d s=4\int_0^{\rho(v)}\frac{\left(f^{-1}(\gamma_s )\right)^{2}}{\left[1+ f^{-1} (\gamma_s)\right]^{1+
  2\kappa}}d s:=4\int_0^{\rho(v)}g(\gamma_s)d s.\end{equation}
Using lemma \ref{lemma}, with probability at least $1-e^{-\tilde{v}}$,

\begin{equation}\label{comp}\int_0^{\tau_{\left(\frac{v}{(1+3\epsilon)c_h}\right)}}g(\gamma_s)d s \leq\frac{\theta_1(v)}{4}\leq \int_0^{\tau_{\left(\frac{v}{(1-3\epsilon) c_h}\right)}}g(\gamma_s)d s.\end{equation}

One can easily check that $g(x)\sim_{\infty}(\kappa
x)^{\frac{1}{\kappa}-2},$ and $g(x)\sim_{-\infty}e^{2x} $. 
In view of this it
is clear that the most important part of the preceding integral will come
from the high values of $\gamma_u.$ To be precise, for $w\in \left[\frac{v}{(1+3\epsilon)c_h},\frac{v}{(1-3\epsilon)c_h}\right]$ and some large constant $A$, we have
\begin{multline}\label{jiunjideux}\int_0^{\tau_w}g(\gamma_s){\bf 1}_{\gamma_s<A}d s=\int_{-\infty}^{A}
g(s)L_{\tau_w}^s d s  \egall w^2 \int_{-\infty}^{A/w}
g(s w)L_{\tau_1}^s d s\\= w^2 \int_{-\infty}^{-A \log(w)^5/w}
g(s w)L_{\tau_1}^s d s+ w^2 \int_{-A\log(w)^5/w}^{A/w}
g(s w)L_{\tau_1}^s d s:=J_1+J_2.\end{multline}

Using statement \ref{mt}, for some constant $C>0$, $\mathbb{P}(J_2>w\log(w)^5)<C e^{-w}.$
Recalling that, under the assumption of theorem $\ref{pas},$ $v\ll \left(\frac{v}{u}\right)^{1/\kappa}$, we get that, for any $\delta>0$, as $t\rightarrow \infty$,
$$\mathbb{P}\left[{J_2}>\delta\left(\frac{v}{u}\right)^{1/\kappa}\right]\leq Ce^{-{\frac{v}{(\log{v})^{10}}}}.$$
We postpone the proof of the following
\begin{lemma}\label{jiun}
for every $\delta>0$, as $t\rightarrow \infty,$
$$\mathbb{P}\left[{J_1}>\delta\left(\frac{v}{u}\right)^{1/\kappa}\right]\leq Ce^{-{\frac{v}{(\log{v})^{10}}}}.$$
\end{lemma}
As a consequence, for every $\delta>0$, as $t\rightarrow \infty,$
\begin{equation}\label{gammaneg}
 \mathbb{P}\left[\int_0^{\tau_w}g(\gamma_s){\bf 1}_{\gamma_s<A}d s>2\delta\left(\frac{v}{u}\right)^{1/\kappa}\right]\leq Ce^{-{\frac{v}{(\log{v})^{10}}}}.
\end{equation}

It remains to deal with $\int_0^{\tau_w}g(\gamma_s){\bf
  1}_{\gamma_s>A}d s.$ 
Due to the equivalent of $g$, for every $\epsilon>0$, for A large enough
\begin{multline}\label{intg}(1-\epsilon)\left(\int_0^{\tau_w}(\gamma_s)^{1/\kappa -2}{\bf
  1}_{\gamma_s>0}d s - I'\right) \leq     \int_0^{\tau_w}g(\gamma_s){\bf
  1}_{\gamma_s>A}d s\\ \leq (1+\epsilon)\int_0^{\tau_w}(\gamma_s)^{1/\kappa -2}{\bf
  1}_{\gamma_s>0}d s,\end{multline}
where $$I':=\int_0^{\tau_w}\gamma_u^{1/\kappa-2}{\bf
  1}_{\gamma_u<A}d u\egall w^{1/{\kappa}}\int_0^{A/w}y^{1/\kappa-2} L^y_{\tau(1)}d y$$
by the same computations as above. Using statement \ref{mt}, for some constant $C'>0$, with probability at least $(1-e^{-C' v})$, $L^y_{\tau(1)}$ is lesser than, say, $100$ on $[0,A/w]$. Therefore
 \begin{equation}\label{iprime} I'\leq 100 w^{1/{\kappa}}\int_0^{A/w}y^{1/\kappa-2}dy<1000 A^{1/\kappa -1} w.\end{equation}

By the same proof as on page 218 of \cite{Jeulin:1981}, the process $$U_s=\int_0^{\tau_s}(\gamma_u)^{1/\kappa -2}{\bf
  1}_{\gamma_s>0}d u$$ is an asymmetric $\kappa$-stable subordinator, more precisely
$$\mathbb{E}\left[\exp{-\frac{\lambda}{2} U_s}\right]=\exp{\left(-s c_\kappa
  \lambda^\kappa\right)},$$
where $c_\kappa=\frac{\pi}{2\kappa \sin{(\pi
  \kappa)}}\left(\frac{\kappa^\kappa}{\Gamma(\kappa)}\right)^2$.
From a result of de Bruijn (see p 221 of \cite{Bertoin:1996}), there
  exists a constant $C_0$ such that
\begin{equation}\label{equivalent 1}\log{\mathbb{P}\left[\frac{U_s}{s^{1/\kappa}}<\left(\frac{1}{u}\right)^{1/\kappa}\right]}=\log{\mathbb{P}\left[U_1<\left(\frac{1}{u}\right)^{1/\kappa}\right]
}    \sim_{\infty}-C_0 u^{\frac{1}{1-\kappa}}.\end{equation}
Similarly, by standard estimates on stable laws, for
$u\rightarrow \infty$, there exists a constant $C_0'$ such that
\begin{equation}\label{equivalent2}\mathbb{P}\left[\frac{U_s}{s^{1/\kappa}}>
  u^{1/\kappa}\right]\sim_{\infty} \frac{C_0'}{u}.\end{equation}
This, together with (\ref{comp}), (\ref{gammaneg}), (\ref{intg}) and (\ref{iprime}), implies that,
for $u\rightarrow \infty$, $u\ll v^{1-\kappa}$ there exists
  positive constants $C_1$ and $C_2$
  such that,
$$\lim_{t\rightarrow\infty}\frac{-\log\mathbb{P}\left[\theta_1(v)<\left(\frac{v}{u}\right)^{1/\kappa}\right]}{u^{\frac{1}{1-\kappa}}}=C_1$$
and for $u\ll e^{v}$, 
$$ \lim_{t\rightarrow \infty}u\mathbb{P}\left[\theta_1(v)>(v u)^{1/\kappa}
  \right]=C_2,$$
where
  $$ C_1=4^{\frac{\kappa}{1-\kappa}}\frac{C_0}{c_h^{\frac{1}{1-\kappa}}}$$ and
  $$C_2=4^{\kappa}\frac{C_0'}{c_h}.$$

\subsubsection{The negative part.}

To finish the proof of Theorem \ref{pas}, we need to deal with
$\theta_2$.
Note that for $\varepsilon>0$,
\begin{equation}\mathbb{P}\left[H_v<\left(\frac{v}{u}\right)^{1/\kappa}\right] \leq
\mathbb{P}\left[\theta_1(v)<\left(\frac{v}{u}\right)^{1/\kappa}\right],\end{equation}
hence the lower bound in (\ref{htasu}) is direct. 

We now turn to the upper bound. We recall that $u$ and $v$ are two functions of $t$ such that $u\ll v^{1-\kappa-\epsilon}.$ This implies in particular that $u\ll v$.  Note that
\begin{equation}\label{hv}\mathbb{P}\left[\theta_1(v)<(1-\varepsilon)\left(\frac{v}{u}\right)^{1/\kappa},\theta_2(v)<\varepsilon\left(\frac{v}{u}\right)^{1/\kappa}\right]\leq\mathbb{P}\left[H(v)<\left(\frac{v}{u}\right)^{1/\kappa}\right].\end{equation}
Using statement \ref{rep}, we obtain
\begin{multline}\label{thetaun}
\mathbb{P}\left[\theta_1(v)<(1-\varepsilon)\left(\frac{v}{u}\right)^{1/\kappa},\theta_2(v)<\varepsilon\left(\frac{v}{u}\right)^{1/\kappa}\right]\\=\mathbb{P}\left[\Upsilon_{2-2\kappa}\left(e^{\Xi_{\kappa}( v)/2}\leadsto 1\right)<\varepsilon\left(\frac{v}{u}\right)^{1/\kappa},\theta_1(v)<(1-\varepsilon)\left(\frac{v}{u}\right)^{1/\kappa}\right].
\end{multline}
By a scaling argument, we get, for $a\geq 1$
\begin{equation}\label{scaling}\mathbb{P}\left(\Upsilon_{2-2\kappa}\left(\sqrt{a}\leadsto 1\right)<a\right)=\mathbb{P}\left(\Upsilon_{2-2\kappa}\left(1\leadsto \frac{1}{\sqrt{a}}\right)<1\right)\geq C>0.\end{equation}
We recall from section \ref{tpp} the representation 
$$e^{\Xi_{\kappa}(v)}-1=f^{-1}(\gamma(\rho(t))).$$
Let $0<\epsilon<\varepsilon/1000$, and $\delta<\varepsilon/3$
we call $A$ the event that the condition of lemma \ref{lemma} is fullfilled, that is
$$A=\left\{\tau_{v/(1+3\epsilon)c_h}<\rho(v)<\tau_{v/(1-3\epsilon)c_h}\right\},$$

Set $\epsilon'\leq (\varepsilon \kappa)^{1/\kappa}/2,$ we introduce the event 
$$B:=\left\{\sup_{\tau_{v/(1+3\epsilon)c_h}<s<\tau_{v/(1-3\epsilon)c_h}} \gamma(s)< \epsilon' \frac{v}{u}\right\}.$$

 Formula 4.1.2 page 185 of \cite{Borodin:2002} (and the Markov property) implies
$$ \mathbb{P}\left[B\right]\geq e^{-\varepsilon'u}$$
for some positive $\varepsilon'$.
We recall from section \ref{tpp} the representation 
$$e^{\Xi_{\kappa}(v)}-1=f^{-1}(\gamma(\rho(t))),$$ where $f^{-1}$ is an increasing function such that $f^{-1}(z)\sim_{\infty}z^{\kappa}/\kappa$.
Therefore for $t$ large enough, on $B\cap A,$  $$e^{\Xi_{\kappa}(v)}<\varepsilon\left(\frac{v}{u}\right)^{1/\kappa}.$$
Recalling equation (\ref{hv}), (\ref{thetaun}), and Lemma \ref{lemma}, we get for $t$ large enough
\begin{multline*}
 \mathbb{P}\left[H(v)<\left(\frac{v}{u}\right)^{1/\kappa}\right]\\\geq \mathbb{P}(B)\mathbb{P}\left[\Upsilon_{2-2\kappa}\left(e^{\Xi_{\kappa}( v)/2}\leadsto 1\right)<\varepsilon\left(\frac{v}{u}\right)^{1/\kappa},\theta_1(v)<(1-\varepsilon)\left(\frac{v}{u}\right)^{1/\kappa}|B\right]\\ \geq 
 \mathbb{P}(B)\mathbb{P}\left[\Upsilon_{2-2\kappa}\left(\sqrt{\varepsilon\left(\frac{v}{u}\right)^{1/\kappa}}\leadsto 1\right)<\varepsilon\left(\frac{v}{u}\right)^{1/\kappa},\right.\\ \left.\theta_1(v)<(1-\varepsilon)\left(\frac{v}{u}\right)^{1/\kappa}|B\right]-\mathbb{P}(B)\mathbb{P}(A^c|B). \end{multline*}
Recalling lemma \ref{lemma} we get $$\mathbb{P}(B)\mathbb{P}(A^c|B)<e^{-\frac{v}{\log v}}.$$
On the other hand, $\Upsilon_{2-2\kappa}\left(\sqrt{\varepsilon\left(\frac{v}{u}\right)^{1/\kappa}}\leadsto 1\right)$ is independent of $B$ and $\theta_1$, and  
$$\mathbb{P}\left[\Upsilon_{2-2\kappa}\left(\sqrt{\varepsilon\left(\frac{v}{u}\right)^{1/\kappa}}\leadsto 1\right)<\varepsilon\left(\frac{v}{u}\right)^{1/\kappa}\right]>C$$
by (\ref{scaling});
therefore the upper bound in (\ref{htasu}) will follow as soon as we show that
$$\lim_{t\rightarrow\infty}\frac{-\log{\mathbb{P}\left[\theta_1(v)<(1-\varepsilon)\left(\frac{v}{u}\right)^{1/\kappa}|B\right]}}{u^{\frac{1}{1-\kappa}}}\leq C_1+\mu(\varepsilon),$$
where $\mu(\varepsilon)\rightarrow 0$ as $\varepsilon \rightarrow 0$.
We recall from equation (\ref{defg}) that
$$g(x)=\frac{\left(f^{-1}(\gamma_s )\right)^{2}}{\left[1+ f^{-1} (\gamma_s)\right]^{1+
  2\kappa}},$$
where $f$ has been defined in (\ref{deff}).

We now recall from equation (\ref{comp}) that, on $A$
$$ \frac{\theta_1(v)}{4}\leq \int_0^{\tau_{\left(\frac{v}{(1-3\epsilon)c_h}\right)}}g(\gamma_s)\mathds{1}_{\gamma_s>0} d s,$$ 
therefore
\begin{multline*}\mathbb{P}\left[\theta_1(v)<(1-\varepsilon)\left(\frac{v}{u}\right)^{1/\kappa}|B\right]\\\geq \mathbb{P}\left[\int_0^{\tau_{\left(\frac{v}{(1-3\epsilon)c_h}\right)}}g(\gamma_s)d s<(1-\varepsilon)\left(\frac{v}{u}\right)^{1/\kappa}|B\right]\\-\mathbb{P}[A^c|B].\\\end{multline*}
 Once again, $ \mathbb{P}[A^c|B]$ is easily bounded.
On the other hand, by Ito's brownian excursion theory (see for example chapter XII of \cite{Revuz:1994}), for every $l\in\mathbb{R},$ $\gamma(\tau_l+t)$ is a brownian motion started at $0$, independent of $\left(\gamma(t)\right)_{t\leq \tau_l}$. Therefore
\begin{multline*}\mathbb{P}\left[\int_0^{\tau_{\left(\frac{v}{(1-3\epsilon)c_h}\right)}}g(\gamma_s)d s<(1-\varepsilon)\left(\frac{v}{u}\right)^{1/\kappa}|B\right]\\ \geq\mathbb{P}\left[\int_0^{\tau_{\left(\frac{v}{(1+3\epsilon)c_h}\right)}}g(\gamma_s)d s<(1-2\varepsilon)\left(\frac{v}{u}\right)^{1/\kappa}|B\right]\\ \mathbb{P}\left[\int_{\tau_{\left(\frac{v}{(1+3\epsilon)c_h}\right)}}^{\tau_{\left(\frac{v}{(1-3\epsilon)c_h}\right)}}g(\gamma_s)d s<\varepsilon\left(\frac{v}{u}\right)^{1/\kappa}|B\right].\end{multline*}

The event in the first probability on the right hand side is independent from $B$, therefore the conditionnal expectation is equal to the expectation and we can apply the results of section \ref{tpp} to get
\begin{equation*}-\log\mathbb{P}\left[\int_0^{\tau_{\left(\frac{v}{(1+3\epsilon)c_h}\right)}}g(\gamma_s)d s <(1-2\varepsilon)\left(\frac{v}{u}\right)^{1/\kappa}|B\right]\leq (C_1 +\mu_1(\varepsilon)+o(1))u^{\frac{1}{1-\kappa}}.\end{equation*}
On the other hand, using the Markov property,
\begin{multline*}\mathbb{P}\left[\int_{\tau_{\left(\frac{v}{(1+3\epsilon)c_h}\right)}}^{\tau_{\left(\frac{v}{(1-3\epsilon)c_h}\right)}}g(\gamma_s)d s<\varepsilon\left(\frac{v}{u}\right)^{1/\kappa}|B\right]\\= \mathbb{P}\left[\int_{0}^{\tau_{\delta v}}g(\gamma_s)d s<\varepsilon\left(\frac{v}{u}\right)^{1/\kappa}|\sup_{0<t<\tau_{\delta v}} \gamma_s<\epsilon' \frac{v}{u}\right],\end{multline*}
where $\delta=\frac{1}{((1-3\epsilon)c_h)}-\frac{1}{((1-3\epsilon)c_h)}$. 
Note that, as the positive and negative excursions are independent, 
$\int_{0}^{\tau_{\delta v}}g^-(\gamma_s)d s$ and B are independent, therefore we only need to bound
\begin{multline*}
\mathbb{P}\left[\int_{0}^{\tau_{\delta v}}g(\gamma_s)d s<\frac{\varepsilon}{2}\left(\frac{v}{u}\right)^{1/\kappa}|\sup_{0<t<\tau_{\delta v}} \gamma_s<\epsilon' \frac{v}{u}\right]\\= \mathbb{P}\left[\int_{0}^{\infty}g(x)L_{\tau_{\delta v}}^x d x<\frac{\varepsilon}{2}\left(\frac{v}{u}\right)^{1/\kappa}|L_{\tau_{\delta v}}^{\alpha}=0\right].
\end{multline*}
 where $\alpha=\epsilon' \frac{v}{u}$.

Intuitively, it seems clear that $\int_{0}^{\infty}g(x)L_{\tau_{\delta v}}^x d x$ will have better chances to be small if $L_{\tau_{\delta v}}^{\alpha}=0$, we are going to give a rigorous proof of that.

Note that, using the second Ray-Knight theorem (Statement \ref{rn}), $L_{\tau_{\delta v}}^x$ is a squared Bessel process of dimension $0$ starting from $\delta v$.
On the other hand, under $\mathbb{P}[\cdot|L_{\tau_{\delta v}}^{\alpha}=0],$ 
$L_{{\tau_{\delta v}}}^x$ is a squared Bessel bridge of dimension $0$ between $\delta v$ and $0$ over time $\alpha$ (we refer to section XI of \cite{Revuz:1994} for the definition and properties of the Bessel bridge).

We are going to use Girsanov's theorem in order to compute the equation solved by the squared Bessel bridge of dimension $0$. 
Let ${\bf P}_x$ and ${\bf P}_{x,0}^{\alpha}$ be respectively the distributions of the Bessel process of dimension $0$ started at $x$ and the distribution of the Bessel bridge of dimension $0$ between $x$ and $0$ over time $\alpha$. Let ${\bf E}_x$ and ${\bf E}_{x,0}^{\alpha}$ be the associated expectations. Let $X_t$ be the canonical process and $\mathcal{F}_t$ its canonical filtration.

Using the Markov property, we get, for every $\mathcal{F}_t$-measurable function $F$, 
\begin{multline*} {\bf  E}_{x,0}^{\alpha}[F(X_s,s\leq t)]=\frac{{\bf  E}_{x}[F(X_s,s\leq t),X_{\alpha}=0]}{{\bf  P}_{x}[X_{\alpha}=0]}\\={\bf  E}_{x}\left[F(X_s,s\leq t)\frac{{\bf  P}_{X_t}[X_{(\alpha-t)}=0]}{{\bf  P}_{x}[X_{\alpha}=0]}\right]:={\bf  E}_{x}\left[F(X_s,s\leq t)h(X_t,t)\right];
 \end{multline*}
where $h(s,t)$ can be explicitely computed (see for example Corollary XI.1.4 of \cite{Revuz:1994}).
We get  
$$h(X_t,t)=\exp{\left(\frac{x}{2\alpha}-\frac{X_t}{2(\alpha-t)}\right)}.$$
Using Ito's Formula, we can transform this expression to get
$$h(X_t,t)=\exp{\left(-\int_{0}^t\frac{1}{2(\alpha-s)}dX_s+\int_0^{t}\frac{X_s}{2(\alpha-s)^2}d s\right)}.$$
Recalling that, under ${\bf P}_x$, $X_t$ is a solution to 
$$dX_t=2\sqrt{X_t}d\beta_t,$$ where $\beta$ is a Brownian motion, we get
$$h(X_t,t)=\exp\left(-\int_{0}^t\frac{\sqrt{X_s}}{(\alpha-s)}d\beta_s+\int_0^{t}\frac{X_s}{2(\alpha-s)^2}d s\right).$$
Therefore, thanks to Girsanov's theorem (see for example Theorem VIII.1.7 of \cite{Revuz:1994}), under ${\bf P}_{x,0}^{\alpha}$,
$$X_t=x+\int_0^t\sqrt{X_s}d\beta_s-2\int_0^{t}\frac{X_s}{(\alpha-s)}d s.$$
Coming back to our original problem, we obtain that, under $\mathbb{P}(\cdot|L_{\tau_{\delta v}}^{\alpha}=0)$, $L_{\tau_{\delta v}}^x$ is a solution to 
$$X_t=\delta v+\int_0^t\sqrt{X_s}d\beta_s-2\int_0^{t}\frac{X_s}{(\alpha-s)}d s.$$
while, under $\mathbb{P}$, $L_{\tau_{\delta v}}^x$ is a solution to 
$$X_t=\delta v+\int_0^t\sqrt{X_s}d\beta_s.$$
Therefore, as there is pathwise uniqueness for these equation (see for example Theorem IX.3.5 of \cite{Revuz:1994}), the comparison theorem (see \cite{watanabe1977comparison}) allows us to construct a couple $(X^{(1)},X^{(2)})$ such that $X^{(1)}$ follows the same distribution as $L_{\tau_{\delta v}}^x$ under $\mathbb{P}$, $X^{(2)}$ follows the same distribution as $L_{\tau_{\delta v}}^x$ under $\mathbb{P}(\cdot|L_{\tau_{\delta v}}^{\alpha}=0)$ and $X^{(1)}\geq X^{(2)}$ almost surely. Then one gets easily that the distribution of $\int_{0}^{\infty}g(x)L_{\tau_{\delta v}}^x d x$ under 
$\mathbb{P}(\cdot|L_{\tau_{\delta v}}^{\alpha}=0)$ is dominated by its distribution under $\mathbb{P}$. Then the upper bound in (\ref{htasu}) follows easily by the results of section \ref{tpp}.\\
 
We now turn to the proof of (\ref{htasd}). We have the trivial inequality
\begin{multline*}\mathbb{P}[\theta_1(v)>(v u)^{1/\kappa}]\leq\mathbb{P}[\theta_1(v)+\theta_2(v)>(v u)^{1/\kappa}] \\
 \leq
\mathbb{P}[\theta_1(v)>(1-\varepsilon)(v u)^{1/\kappa}]+ \mathbb{P}[\theta_2(v)>\varepsilon(v u)^{1/\kappa}],\end{multline*}
therefore the lower bound is direct. To get the upper bound, note that $\theta_2(v)$ is increasing, so we have to show that for every $\varepsilon>0,$ and some $s>0$,
$$\mathbb{P}\left[\Upsilon_{2-2\kappa}\left(e^{\Xi_{\kappa}(v+s)/2}\leadsto 1\right)>\varepsilon\left({v}{u}\right)^{1/\kappa}\right]=o\left(\frac{1}{u}\right).$$
Recalling the diffusion $Z_t$ from the last part, we need to bound
\begin{multline}\label{int}
\mathbb{P}\left[\Upsilon_{2-2\kappa}\left(\sqrt{Z_{v+s}+1}\leadsto 1\right)>\varepsilon\left({v}{u}\right)^{1/\kappa}\right]
\\=\int_0^{\infty}\mathbb{P}\left[\Upsilon_{2-2\kappa}\left(\sqrt{z+1}\leadsto 1\right)>\varepsilon\left({v}{u}\right)^{1/\kappa}\right]d\mu_{v+s} (z),
\end{multline}
where $\mu_v(y)$ is the distribution of $Z_v$.
By scaling, 
\begin{multline*}\mathbb{P}\left[\Upsilon_{2-2\kappa}\left(\sqrt{z+1}\leadsto 1\right)>\varepsilon\left({v}{u}\right)^{1/\kappa}\right]
\\= \mathbb{P}\left[\Upsilon_{2-2\kappa}\left(1\leadsto \frac{1}{\sqrt{z+1}}\right)>\varepsilon\frac{(v u)^{1/\kappa}}{z+1}\right]\leq\mathbb{P}\left[\Upsilon_{2-2\kappa}\left(1\leadsto 0\right)>\varepsilon\frac{(v u)^{1/\kappa}}{z+1}\right]. \end{multline*}
It is known (see for example \cite{Werner:1995}
page 40) that $\Upsilon_{2-2\kappa}\left(1\leadsto 0\right)$ has the same distribution as $\frac{1}{2\Gamma}$ where $\Gamma$ follows a distribution ${\Gamma}(\kappa,1)$, therefore, easy computations leads to 
$$\mathbb{P}\left[\Upsilon_{2-2\kappa}\left(1\leadsto 0\right)>\varepsilon\frac{(v u)^{1/\kappa}}{z+1}\right]\leq \left(\frac{1}{\kappa\Gamma(\kappa)(2\epsilon)^{\kappa}}\frac{(1+z)^{\kappa}}{u v}\right)\wedge1.$$
Recalling (\ref{int}), we have, for all $A>0$
\begin{equation*}\mathbb{P}\left(\theta_2>\varepsilon(uv)^{1/\kappa}\right)\leq \int_0^A \frac{1}{\kappa\Gamma(\kappa)(2\epsilon)^{\kappa}}\frac{(1+z)^{\kappa}}{uv}d\mu_{v+s}(z)+\int_A^\infty 
d\mu_{v+s}(z)\end{equation*}

Using for example exercise VII.3.20 of \cite{Revuz:1994}, the diffusion $Z_t$ has speed measure $d m(z)=\frac{2}{(1+z)^{1+\kappa}}d z$, so by Theorem 54.4 of \cite{Rogers:1986} and a change of variable in order to lift the natural scale assumption, for any $\phi$ bounded and measurable, 
$$\int_0^{\infty}\phi(z)d\mu_{v+s}(z)\rightarrow_{s\rightarrow \infty}\int_0^{\infty} \phi(z)\pi(d z),$$
with $\pi(d z) =\frac{m(d z)}{2\kappa}$. Therefore as $s$ goes to infinity, and for some finite constant $c(\varepsilon),$ 
$$
\mathbb{P}\left(\theta_2>\varepsilon(u v)^{1/\kappa}\right)\leq\frac{c(\varepsilon)}{u v}\log(1+A)+(1+A)^{-\kappa}.
$$
Now, taking $A$ such that $(1+A)\gg u^{1/\kappa}$ and $\log(1+A)\ll v$ (this is possible due to the assumptions on $u$ and $v$), we get the upper bound in (\ref{htasd}).

\subsection{Proof of Theorem \ref{as}.}
In this section we use the results for the hitting times to get the results for the diffusion itself.
We begin with the proof of (\ref{ase1}).
We have the trivial inequality
$$\mathbb{P}\left(X_t>t^\kappa u \right)\leq
\mathbb{P}\left[H(t^\kappa u)<t \right]; $$
by taking $v=t^\kappa u$ in Theorem \ref{pas}, we get the upper bound in (\ref{ase1}). The condition
$u\ll v^{1-\kappa}$ becomes $u\ll t^{1-\kappa}$.

To get the lower bound, note that, for every $\varepsilon>0,$
\begin{multline}\label{le}\log \mathbb{P}\left(X_t>t^\kappa u \right) \\\geq
  \log \left[\mathbb{P}\left[H((1+\varepsilon)t^\kappa u)<t
    \right]\mathbb{P}\left(X_t>t^\kappa u |H((1+\varepsilon)t^\kappa
  u)<t\right) \right]\\
\geq -C_1((1+\varepsilon)u)^{\frac{1}{1-\kappa}} + \log\mathbb{P}\left(X_t>t^\kappa u |H((1+\varepsilon)t^\kappa
  u)<t\right).
\end{multline}
The bound in the first term coming from (\ref{htasu}).
To treat the second term, note that
\begin{multline*}
\mathbb{P}\left(X_t<t^\kappa u |H((1+\varepsilon)t^\kappa
  u)<t\right)\leq E\left[P_W^{(1+\varepsilon)t^\kappa
  u}\left[\inf_{s>0} X_s < t^\kappa
  u\right]\right]=\\\mathbb{P}\left[\inf_{s>0} X_s <- \varepsilon t^\kappa
  u\right],\end{multline*}
by invariance of the environment.

By \cite{Kawazu:1993},
$$\mathbb{P}\left[\inf_{t>0} X_t <- u \right]\leq C x^{-3/2}\exp{\left(-(\kappa/2)^2 x/2\right)},$$ (note that $c$ in K. Kawazu and
  H. Tanaka's article corresponds to $-\kappa/2$ in our
  setting). Therefore we get easily that, for $t$ large enough, 
\begin{equation}\label{xt}\mathbb{P}\left(X_t<t^\kappa u |H((1+\varepsilon)t^\kappa
  u)<t\right)<1/2.
\end{equation}
The lower bound in (\ref{ase1}) then follows from equations (\ref{le}) and (\ref{xt}).

To prove (\ref{ase2}), we use the fact that, for every $\varepsilon>0$,
\begin{multline*}\mathbb{P}\left[H\left(\frac{t^\kappa}{u}>t\right)\right]\leq\mathbb{P}\left[X_t<\frac{t^\kappa}{u}\right]\\
\leq\mathbb{P}\left[H\left(\frac{(1+\varepsilon)t^\kappa}{u}\right)> t\right]+\mathbb{P}\left[X_t<\frac{t^\kappa}{u};H\left(\frac{(1+\varepsilon)t^\kappa}{u}\right)< t\right].
\end{multline*}
Taking $v=t^\kappa/u$, Theorem \ref{pas} implies the lower
bound, and the upper bound follows easily by the same argument as
before.

It remains to prove Lemma \ref{lemma} and Lemma \ref{jiun}
\subsection{Proof of Lemmas \ref{lemma} and \ref{jiun}.}
We begin with the proof of Lemma \ref{lemma}. It will turn out that once the tools for this Lemma will we introduced, Lemma \ref{jiun} will be quite obvious.
We recall from equation (\ref{defalpha}) that 
$$\alpha=\rho(t)^{-1}=\int_0^t h(\gamma_s)d s,$$
where $h$ is some positive, integrable function. We have
$$\alpha_{\tau_t}=\int_{-\infty}^{\infty} h(x)L_{\tau_t}^x d x=t c_h +
\int_{-\infty}^{\infty} h(x)\left(L_{\tau_t}^x-t\right)d x.$$
Our result then follows as soon as we show that, for $t$ large enough
$$\mathbb{P}\left[\left|\int_{-\infty}^{\infty}
  h(x)\left(L_{\tau_t}^x-t\right)d x>3t\epsilon\right|\right]<\exp\left(-w\right).$$
Let $s$ such that $s\rightarrow \infty$ and $t/s^4 \gg w,$ then
\begin{multline*}
\int_{-\infty}^{\infty}
  h(x)\left(L_{\tau_t}^x-t\right)d x=\int_{-s}^{s}
  h(x)\left(L_{\tau_t}^x-t\right)d x\\
+\int_{s}^{\infty}
  h(x)\left(L_{\tau_t}^x-t\right)d x\\
+\int_{-\infty}^{s}
  h(x)\left(L_{\tau_t}^x-t\right)d x\\:=I_1 + I_2 + I_3.\\
\end{multline*}

By a scaling argument, and using the fact that $h$ is bounded, we have 
$$|I_1|\leq C\int_{-s}^{s}\left|L_{\tau_t}^x-t\right|d x\egall
C t\int_{-s}^{s}\left|L_{\tau_1}^{x/t}-1\right|d x
=C t^2\int_{-s/t}^{s/t}\left|L_{\tau_1}^{y}-1\right|d y.$$
Then, for $t$ large enough,
\begin{multline*}\mathbb{P}(|I_1|>t\epsilon)\leq\mathbb{P}\left(\sup_{y\in[-s/t,s/t]}|L_{\tau_1}^{y}-1|>\frac{\epsilon}{2Cs}\right)\\\leq 2\mathbb{P}\left(\sup_{y\in[0,s/t]}|L_{\tau_1}^{y}-1|>\frac{\epsilon}{2Cs}\right),\end{multline*}
the last bound coming from the symmetry of $L_{\tau_1}^y$ in $y$.
 On the other hand, using statement \ref{rn}, $L_{\tau_1}^{y}$ is a squared Bessel process of dimension $0$ started from $1$, therefore using statement \ref{mt} with $\delta=\frac{\epsilon}{2Cs},$ $v=s/t$, we get
\begin{equation*}\mathbb{P}(|I_1|>t\epsilon) \leq C's\exp\left(-\frac{\epsilon^2t}{s^3}\right)\leq \exp(-w).\end{equation*}
 It is clear that, for large $t$, $\mathbb{P}(|I_3| \geq
 t\epsilon)\leq\mathbb{P}(|I_2| \geq t\epsilon).$
To bound $I_3$, we note that, for $t$ large enough,
$$|I_3|\leq 2\left(\int_{-\infty}^{s}
  \frac{L_{\tau_t}^x}{x^2}d x+t\int_{-\infty}^{s}\frac{1}{x^2}d x\right)\egall 2\left(\frac{t}{s}+\int_{s/t}^{\infty}
  \frac{L_{\tau_1}^x}{x^2}d x\right), $$
by the same scaling argument. The first part is negligible, and, using statement \ref{rn},
$$\int_{s/t}^{\infty}
  \frac{L_{\tau_1}^x}{x^2}d x=\int_{s/t}^{\infty}
  \frac{Z_x}{x^2}d x,$$ where $Z_t$ is a squared Bessel process of dimension $0$ started at 1.  The following result from
  \cite{Pitman:1982} allows us to compute the Laplace transform of
  this random variable.
\begin{statement}[J. Pitman and M.Yor]\label{pit}
Let $Z_t$ be a squared Bessel process of dimension $d$, starting from
$x$, and $\mu$ a positive (Radon) measure on $(0,\infty)$ such that,
for all $n$, $\mu(0,n)<\infty$. Then one has
$$E\left[\exp{\left(-\int Z_t
    d\mu(t)\right)}\right]=\phi_\mu(\infty)^{d/2}
    \exp{\left(\frac{x}{2} \phi_\mu'(0)\right)},$$
where  $\phi_\mu$ is the unique decreasing and convex solution of
$$\frac{1}{2} \phi'' =\mu.\phi\;\text{on } (0,\infty),\,\phi(0)=1.$$
\end{statement}
We note $\eta = s/t$, and $A_t=\int_{\eta}^{\infty}
  \frac{L_{\tau_1}^x}{x^2}d x$.
The preceding statement implies that
$$\mathbb{E}\left[\exp{-\lambda A_t}\right]=\exp{\left(\frac{1}{2} \phi_\mu'(0)\right)},$$
 where $\phi_\mu$ is the solution of:
$$\phi''(x)=2 \lambda \frac{\phi(x)}{x^2} {\bf 1}_{x\geq \eta}.$$
A decreasing solution on $(\eta,\infty)$ of this equation is
$$\phi(x)=C\left(\frac{x}{\eta}\right)^{\frac{1-\sqrt{1+8\lambda}}{2}}.$$
The condition $\phi(0) = 1$ and the fact that $\phi'$ is constant on
$[0,\eta]$ implies that
$C\left(1-\frac{1-\sqrt{1+8\lambda}}{2}\right)=1$, thus

$$\mathbb{E}\left[\exp{-\lambda
    A_t}\right]=\exp{\left(\frac{1-\sqrt{1+
    8\lambda}}{2(1+\sqrt{1+8\lambda})\eta}\right)},$$
As this function is analytic, for some $\lambda>0$ (not depending on $t$),
$$\mathbb{E}\left[\exp{\lambda
    A_t}\right]=\exp{\left(\frac{1-\sqrt{1-8\lambda}}{2(1+\sqrt{1-8\lambda})\eta}\right)},$$
then 
$$\mathbb{P}(I_2>\epsilon t)\leq
\exp{\left(\frac{1-\sqrt{1-8\lambda}}{2(1+\sqrt{1-8\lambda})\eta}-\lambda
\epsilon t\right)},$$ from which the result follows, as $1/\eta\ll t.$

Let us now prove Lemma \ref{jiun}.
We recall from (\ref{jiunjideux}) that, 
$$J_1=w^2 \int_{-\infty}^{-A\log(w)^5/w}
g(s w)L_{\tau_1}^s d s \leq 2\int_{-\infty}^{-A\log(w)^5/w}
\frac{1}{s^2}L_{\tau_1}^s d s.$$
Then the proof follows easily as a corollary of the proof of Lemma \ref{lemma}.
\section{Quenched slowdown.}\label{quenchedslowdown}
We now turn to the proof of Theorem \ref{pqsl}. As before we first recall some useful facts.
\subsection{Preliminary statements.}
We recall the time change representation of $X_t$ (see, for example \cite{Hu:1999})
$$X_t=A_\kappa^{-1}\left(B(T_\kappa^{-1}(t))\right),$$ where
$$A_\kappa(x)=\int_0^x e^{W\kappa(y)}d y,$$
$$T_\kappa(t)=\int_0^t e^{-2W_\kappa(A_\kappa^{-1}(B(s)))}d s,$$
and $B$ is a standard Brownian motion.

We also need a result about Sturm-Liouville equations.
Let $V(t)$ be a positive function of $t\geq 0$, and $\bar{V}(t)=\int_{0}^t V(u) d u.$ We are interested in the solution  of the differential equation
\begin{equation}\label{sl}z''(t)=-\lambda V(t)z(t),\, t\geq 0,\:\: z(0)=1, \,z'(0)=0.\end{equation}
We have the following statement from \cite{Bobkov:2002} (corollary 3.2)
\begin{statement}\label{bobkov}
Let $\lambda(V)$ be the supremum of all $\lambda>0$ for which a solution to the problem (\ref{sl}) is positive in $[0,1)$, then
$$\sup_{0<t<1}(1-t)\bar{V}(t)\leq \frac{1}{\lambda(V)}\leq 4 \sup_{0<t<1}(1-t)\bar{V}(t) .$$
\end{statement}

We recall the following inequality from lemma 1.1.1 of
\cite{Csorgo:1981}
\begin{statement}\label{Csorgo}
Let $\gamma(t)$ be a one-dimensional brownian motion, then
$$P\left( \sup_{0\leq s_1<s_2<t,
  s_2-s_1<u}|\gamma(s_2)-\gamma(s_1)|>\frac{x}{2}\right)\leq
  c\frac{t}{u} \exp{-\frac{x^2}{9u}}.$$
\end{statement}
We finish with a useful lemma
\begin{lemma}\label{Laplace}
let $a>0$, and $\mu$ a Radon measure on $[0,a]$, and suppose there exists $\phi$ a positive
solution of the
Sturm-Liouville equation
\begin{equation}\label{sl'}\phi''=-\phi\mu,\, t\geq 0,\:\: \phi(a)=1, \,\phi'(a)=0.\end{equation}
Let $X_t$ be a squared Bessel process of dimension $\delta$, starting at x, then
$$E\left[\exp\left(\int_0^a X_t d\mu(t)\right)\right]\leq
\phi(0)^{-\delta/2}\exp{\left(\frac{1}{2}\frac{\phi'(0)}{\phi(0)}x \right)}.$$
\end{lemma}
\textbf{Remark:} This lemma is a extension of Statement \ref{pit}, but we do not get equality in this case.

\textbf{Proof:}
Let $F_\mu(t)=\phi'(t)/\phi(t)$, by the concavity of $\phi$ this is a
right continuous and
decreasing function, thus we can apply the integration by parts
formula to get
$$F_\mu(t)X_t =F_\mu(0) x+\int_0^{t} F_\mu(s) d X_s +\int_0^t X_s
d F_\mu(s).$$
Using (\ref{sl'}), we can compute the last part
\begin{multline*}\int_0^t X_s d F_\mu(s)=\int_0^t X_s \frac {d\phi'(s)}{\phi(s)}
-\int_0^t \frac{\phi'(s)d\phi(s)}{\phi(s)^2}
\\=-\int_0^t X_s d\mu(s)-\int_0^t X_s F_\mu(s)^2d s.
\end{multline*}
Recalling that $M_t=X_t-\delta t$ is a local martingale, we set 
$$
Z_\mu
(t)=\exp\left(\frac{1}{2}\int_0^t
F_\mu(s)d M_s-\frac{1}{2}\int_0^t X_s F_\mu(s)^2 d s\right),
$$
which is a positive local martingale, hence a supermartingale.
Using the previous computation, we get
$$Z_\mu(t)=\exp\left(\frac{1}{2}\left[F_\mu(t)X_t-F_\mu(0)x-\delta\int_0^t
    F_\mu(s)d s + \int_0^t X_s d\mu(s)\right]\right).$$
As $Z_\mu$ is a supermartingale, $E[Z_\mu(a)]\leq E[Z_\mu(0)]=1$. Therefore
the result follows easily.

\subsection{Quenched slowdown for the hitting time.}
In this section we show (\ref{pqsl2}). The idea of the proof is to decompose the environment in valleys of a certain size, then to study the process of the valleys visited and the time spent in the valleys. We first give a formal definition of what a valley is. 
For $t>0$, $v>0$ and $i \in \mathbb{N}$
, we set $K_0=-\lfloor t\rfloor$, and
\begin{multline*}
K_{i+1}=\inf\left\{x>K_{i}, W_{\kappa}(K_{i})-\inf_{y\in[K_i,x]} W_{\kappa}(y)>
  \frac{3}{ \kappa} \log \lfloor t\rfloor,\right.\\\left. W_\kappa(x)\geq\sup_{y>x}
  W_{\kappa}(y)-1\right\}.\end{multline*}
  $K_i$ is finite almost surely, due to the transience of the drifted brownian motion.
The intervals $[K_i,K_{i+1}]$ will be called ``valleys''. An example of such valleys is given in figure 2.
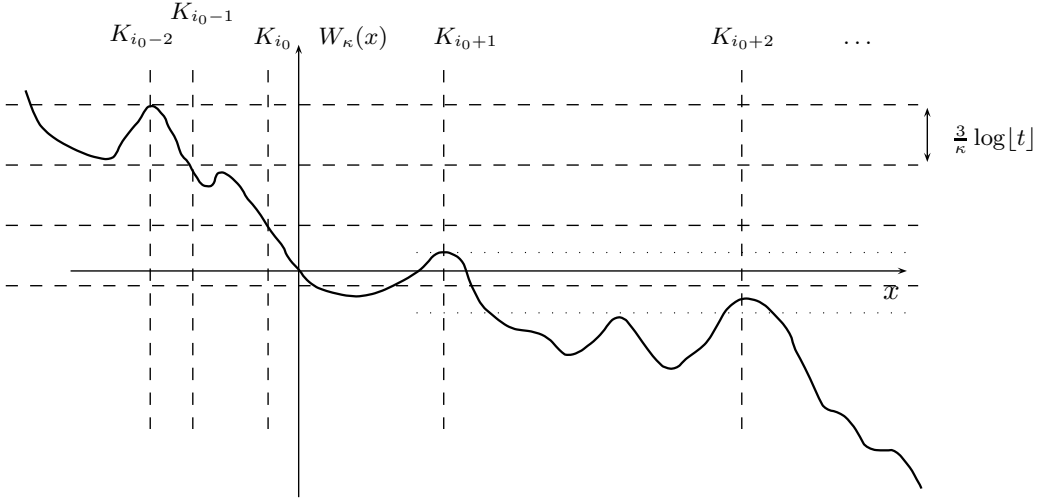
\begin{figure}[h!]\centering
 \scalebox{1} 
{
\begin{pspicture}(0.5,-3.146)(17.2,3.186)
\rput(3.85,-0.146){\psaxes[linewidth=0.018,arrowsize=0.05291667cm 2.0,arrowlength=1.4,arrowinset=0.4,labels=none,ticks=none,ticksize=0.10583333cm]{->}(0,0)(-3,-3)(8,3)}

\pscustom[linewidth=0.03]
{
\newpath
\moveto(0.26,2.246)
\lineto(0.33,2.046)
\curveto(0.365,1.946)(0.43,1.811)(0.46,1.776)
\curveto(0.49,1.741)(0.545,1.686)(0.57,1.666)
\curveto(0.595,1.646)(0.65,1.606)(0.68,1.586)
\curveto(0.71,1.566)(0.81,1.511)(0.88,1.476)
\curveto(0.95,1.441)(1.07,1.391)(1.12,1.376)
\curveto(1.17,1.361)(1.26,1.341)(1.3,1.336)
\curveto(1.34,1.331)(1.4,1.346)(1.42,1.366)
\curveto(1.44,1.386)(1.47,1.436)(1.48,1.466)
\curveto(1.49,1.496)(1.53,1.566)(1.56,1.606)
\curveto(1.59,1.646)(1.64,1.711)(1.66,1.736)
\curveto(1.68,1.761)(1.715,1.801)(1.73,1.816)
\curveto(1.745,1.831)(1.78,1.886)(1.8,1.926)
\curveto(1.82,1.966)(1.855,2.016)(1.87,2.026)
\curveto(1.885,2.036)(1.92,2.041)(1.94,2.036)
\curveto(1.96,2.031)(1.995,2.006)(2.01,1.986)
\curveto(2.025,1.966)(2.055,1.931)(2.07,1.916)
\curveto(2.085,1.901)(2.105,1.866)(2.11,1.846)
\curveto(2.115,1.826)(2.13,1.786)(2.14,1.766)
\curveto(2.15,1.746)(2.18,1.696)(2.2,1.666)
\curveto(2.22,1.636)(2.25,1.571)(2.26,1.536)
\curveto(2.27,1.501)(2.29,1.441)(2.3,1.416)
\curveto(2.31,1.391)(2.34,1.351)(2.36,1.336)
\curveto(2.38,1.321)(2.415,1.271)(2.43,1.236)
\curveto(2.445,1.201)(2.475,1.146)(2.49,1.126)
\curveto(2.505,1.106)(2.53,1.066)(2.54,1.046)
\curveto(2.55,1.026)(2.58,0.996)(2.6,0.986)
\curveto(2.62,0.976)(2.665,0.971)(2.69,0.976)
\curveto(2.715,0.981)(2.745,1.021)(2.75,1.056)
\curveto(2.755,1.091)(2.775,1.136)(2.79,1.146)
\curveto(2.805,1.156)(2.84,1.161)(2.86,1.156)
\curveto(2.88,1.151)(2.93,1.121)(2.96,1.096)
\curveto(2.99,1.071)(3.045,1.016)(3.07,0.986)
\curveto(3.095,0.956)(3.13,0.901)(3.14,0.876)
\curveto(3.15,0.851)(3.18,0.806)(3.2,0.786)
\curveto(3.22,0.766)(3.27,0.706)(3.3,0.666)
\curveto(3.33,0.626)(3.38,0.551)(3.4,0.516)
\curveto(3.42,0.481)(3.465,0.411)(3.49,0.376)
\curveto(3.515,0.341)(3.57,0.271)(3.6,0.236)
\curveto(3.63,0.201)(3.665,0.146)(3.67,0.126)
\curveto(3.675,0.106)(3.69,0.071)(3.7,0.056)
\curveto(3.71,0.041)(3.73,0.011)(3.74,-0.0040)
\curveto(3.75,-0.019)(3.79,-0.064)(3.82,-0.094)
\curveto(3.85,-0.124)(3.89,-0.174)(3.9,-0.194)
\curveto(3.91,-0.214)(3.955,-0.264)(3.99,-0.294)
\curveto(4.025,-0.324)(4.085,-0.364)(4.11,-0.374)
\curveto(4.135,-0.384)(4.19,-0.404)(4.22,-0.414)
\curveto(4.25,-0.424)(4.35,-0.449)(4.42,-0.464)
\curveto(4.49,-0.479)(4.605,-0.489)(4.65,-0.484)
\curveto(4.695,-0.479)(4.78,-0.464)(4.82,-0.454)
\curveto(4.86,-0.444)(4.945,-0.414)(4.99,-0.394)
\curveto(5.035,-0.374)(5.155,-0.314)(5.23,-0.274)
\curveto(5.305,-0.234)(5.42,-0.159)(5.46,-0.124)
\curveto(5.5,-0.089)(5.565,-0.019)(5.59,0.016)
\curveto(5.615,0.051)(5.675,0.091)(5.71,0.096)
\curveto(5.745,0.101)(5.8,0.101)(5.82,0.096)
\curveto(5.84,0.091)(5.895,0.066)(5.93,0.046)
\curveto(5.965,0.026)(6.015,-0.034)(6.03,-0.074)
\curveto(6.045,-0.114)(6.07,-0.189)(6.08,-0.224)
\curveto(6.09,-0.259)(6.12,-0.334)(6.14,-0.374)
\curveto(6.16,-0.414)(6.2,-0.494)(6.22,-0.534)
\curveto(6.24,-0.574)(6.29,-0.639)(6.32,-0.664)
\curveto(6.35,-0.689)(6.415,-0.744)(6.45,-0.774)
\curveto(6.485,-0.804)(6.555,-0.854)(6.59,-0.874)
\curveto(6.625,-0.894)(6.69,-0.919)(6.72,-0.924)
\curveto(6.75,-0.929)(6.825,-0.939)(6.87,-0.944)
\curveto(6.915,-0.949)(7.005,-0.974)(7.05,-0.994)
\curveto(7.095,-1.014)(7.165,-1.064)(7.19,-1.094)
\curveto(7.215,-1.124)(7.26,-1.174)(7.28,-1.194)
\curveto(7.3,-1.214)(7.335,-1.244)(7.35,-1.254)
\curveto(7.365,-1.264)(7.42,-1.264)(7.46,-1.254)
\curveto(7.5,-1.244)(7.57,-1.209)(7.6,-1.184)
\curveto(7.63,-1.159)(7.675,-1.124)(7.69,-1.114)
\curveto(7.705,-1.104)(7.755,-1.059)(7.79,-1.024)
\curveto(7.825,-0.989)(7.875,-0.929)(7.89,-0.904)
\curveto(7.905,-0.879)(7.945,-0.829)(7.97,-0.804)
\curveto(7.995,-0.779)(8.05,-0.759)(8.08,-0.764)
\curveto(8.11,-0.769)(8.155,-0.794)(8.17,-0.814)
\curveto(8.185,-0.834)(8.22,-0.879)(8.24,-0.904)
\curveto(8.26,-0.929)(8.3,-0.984)(8.32,-1.014)
\curveto(8.34,-1.044)(8.4,-1.129)(8.44,-1.184)
\curveto(8.48,-1.239)(8.54,-1.314)(8.56,-1.334)
\curveto(8.58,-1.354)(8.615,-1.389)(8.63,-1.404)
\curveto(8.645,-1.419)(8.7,-1.439)(8.74,-1.444)
\curveto(8.78,-1.449)(8.845,-1.424)(8.87,-1.394)
\curveto(8.895,-1.364)(8.94,-1.319)(8.96,-1.304)
\curveto(8.98,-1.289)(9.035,-1.249)(9.07,-1.224)
\curveto(9.105,-1.199)(9.165,-1.149)(9.19,-1.124)
\curveto(9.215,-1.099)(9.26,-1.039)(9.28,-1.004)
\curveto(9.3,-0.969)(9.34,-0.899)(9.36,-0.864)
\curveto(9.38,-0.829)(9.41,-0.774)(9.42,-0.754)
\curveto(9.43,-0.734)(9.455,-0.694)(9.47,-0.674)
\curveto(9.485,-0.654)(9.525,-0.614)(9.55,-0.594)
\curveto(9.575,-0.574)(9.62,-0.544)(9.64,-0.534)
\curveto(9.66,-0.524)(9.71,-0.514)(9.74,-0.514)
\curveto(9.77,-0.514)(9.84,-0.529)(9.88,-0.544)
\curveto(9.92,-0.559)(9.99,-0.604)(10.02,-0.634)
\curveto(10.05,-0.664)(10.14,-0.759)(10.2,-0.824)
\curveto(10.26,-0.889)(10.33,-1.004)(10.34,-1.054)
\curveto(10.35,-1.104)(10.38,-1.189)(10.4,-1.224)
\curveto(10.42,-1.259)(10.475,-1.364)(10.51,-1.434)
\curveto(10.545,-1.504)(10.61,-1.644)(10.64,-1.714)
\curveto(10.67,-1.784)(10.72,-1.889)(10.74,-1.924)
\curveto(10.76,-1.959)(10.825,-1.999)(10.87,-2.004)
\curveto(10.915,-2.009)(11.0,-2.049)(11.04,-2.084)
\curveto(11.08,-2.119)(11.15,-2.204)(11.18,-2.254)
\curveto(11.21,-2.304)(11.26,-2.394)(11.28,-2.434)
\curveto(11.3,-2.474)(11.385,-2.519)(11.45,-2.524)
\curveto(11.515,-2.529)(11.61,-2.529)(11.64,-2.524)
\curveto(11.67,-2.519)(11.74,-2.564)(11.78,-2.614)
\curveto(11.82,-2.664)(11.875,-2.734)(11.89,-2.754)
\curveto(11.905,-2.774)(11.935,-2.824)(11.95,-2.854)
\curveto(11.965,-2.884)(11.995,-2.939)(12.01,-2.964)
\curveto(12.025,-2.989)(12.04,-3.019)(12.04,-3.034)
}

\rput(3.52,2.936){\footnotesize $K_{i_0}$}

\rput(4.57,2.936){\footnotesize $W_{\kappa}(x)$}

\rput(11.65,-0.439){$x$}
\psline[linewidth=0.018cm,linestyle=dashed,dash=0.16cm 0.16cm](0.0,1.256)(12.0,1.256)
\psline[linewidth=0.018cm,linestyle=dashed,dash=0.16cm 0.16cm](0.0,2.056)(12.0,2.056)
\psline[linewidth=0.018cm,linestyle=dashed,dash=0.16cm 0.16cm](0.0,0.456)(12.0,0.456)
\psline[linewidth=0.018cm,linestyle=dashed,dash=0.16cm 0.16cm](0.0,-0.344)(12,-0.344)
\psline[linewidth=0.018cm,arrowsize=0.05291667cm 2.0,arrowlength=1.4,arrowinset=0.4]{<->}(12.12,2.016)(12.12,1.296)

\rput(13.,1.596){\footnotesize $\frac{3}{\kappa}\log \lfloor t \rfloor$}
\psline[linewidth=0.018cm,linestyle=dashed,dash=0.16cm 0.16cm](1.9,2.516)(1.9,-2.244) 

\psline[linewidth=0.018cm,linestyle=dashed,dash=0.16cm 0.16cm](3.45,2.516)(3.45,-2.244)
\psline[linewidth=0.018cm,linestyle=dashed,dash=0.16cm 0.16cm](5.76,2.516)(5.76,-2.244)
\psline[linewidth=0.018cm,linestyle=dashed,dash=0.16cm 0.16cm](9.68,2.516)(9.68,-2.244)
\psline[linewidth=0.018cm,linestyle=dashed,dash=0.16cm 0.16cm](2.46,2.516)(2.46,-2.244)
\rput(6.05,2.936){\footnotesize $K_{i_0+1}$}

\rput(9.69,2.936){\footnotesize $K_{i_0+2}$}

\rput(2.59,3.266){\footnotesize $K_{i_0-1}$}

\rput(1.79,2.966){\footnotesize $K_{i_0-2}$}
\psline[linewidth=0.018cm,linestyle=dotted,dotsep=0.16cm](5.4,0.096)(12.0,0.096)
\psline[linewidth=0.018cm,linestyle=dotted,dotsep=0.16cm](5.4,-0.704)(12.0,-0.704)
\rput(11.23,2.881){\footnotesize \dots} 
\end{pspicture}
 } \label{figure2}
\caption{Decomposition in Valleys}
\end{figure}

We introduce the sequence defined, for $k\geq 0$ by
\begin{align*}
s_0&=0\\
s_{k+1}&=\inf\{t>s_k, X_t \in \{K_j,j\geq 0\}\}.
\end{align*}
We call $Y_k=X_{s_k}$, ${l_t}=\max\{i:s_i<H(v)\}$ and 
$$\xi(i)=\sharp\left\{j\in
  [0,l_t],Y_j=K_{i+1},Y_{j+1}=K_i\right\}.$$
We set $i_0=\max\{j, K_j< 0\}$ and $i_1=\max\{j, K_j< v \}$. By convention we note $K_{i_1+1}=v$. 
Let
\begin{equation}\label{defb}\mathcal{B}=\sum_{i=1}^{i_1-1} \xi(i)\end{equation} denote the number of times the "walk" $Y_k$ backtracks.
 Let
  $\theta(t)$ be the time-shift associated to the diffusion, we
  set for $0\leq i < i_1$
$$next(i)=\inf\{t\geq 0: X_t=K_i, H({K_{i+1}})\circ
  \theta(t)<H(K_{i-1})\circ\theta(t)\}$$ and
 $$H^{next}(i)=H(K_{i+1})\circ \theta(next(i))-next(i).$$
We have the following decomposition of ${H_v}:$
$$H(v)=H_{init}+H_{dir}+H_{back}+H_{left}+H_{right},$$
where
\begin{equation}
H_{init}=\left\{\begin{array}{c c}H(K_{i_0+1}) & \text{if
}H(K_{i_0+1})<H(K_{i_0})\\
H(K_{i_0})+H^{next}(i_0)\circ\theta(H(K_{i_0})) & \text{else}
    \end{array} \right.,
\label{hinit}\end{equation}
is the time the diffusion takes to get to $K_{i_0+1},$
\begin{equation}\label{hleft}
H_{left}=\int_0^t {\bf 1}_{X_t<K_1} d t ,
\end{equation}
is the time the diffusion spends at the left of $K_1$,
\begin{equation}\label{hright}H_{right}=H(v)\circ \theta(next(i_1)) - next(i_1),\end{equation}
 is the time spent to get from $K_{i_1}$ to $v$
\begin{equation}\label{hdir} H_{dir}= \sum_{i=i_0+1}^{i_1-1} H^{next}(i)\end{equation} 
is the time used for the direct crossings of the valleys and
\begin{multline}\label{hback}H_{back}=\sum_{i=i_0+1}^{i_1-1} \sum_{j=0}^{l_t} {\bf
    1}_{Y_j=K_{i+1},Y_{j+1}=K_j}\\
\times \left( H(K_i)\circ \theta(s_j) - s_j + H^{next}(i)
\circ\theta(H({K_i})\circ \theta (s_j))\right)
\end{multline}
is the time ``lost'' as a consequence of the different backtracks of $Y_k$.

We introduce $D_i=\sup_{K_i<s<t<K_{i+1}} W_{\kappa}(t)-W_{\kappa}(s)$, to which we will refer as the ``depth'' of the valley $[K_i,K_{i+1}]$,
and $$N(s,t)=\{ i\geq 1, [K_i,K_i+1)\cap [s,t)\neq \emptyset \}.$$ Note that, as seen on figure 2 there are some valleys of depth $0$. 

We have the following lemmas, whose proof will be postponed
\begin{lemma}[environment estimates]\label{enve}
Let $v=t^\nu$ and $\epsilon>0$.
${\cal P}$-almost surely, for $m>m_0$, for $t$ large enough, $W\in \Omega$
where $\Omega=\Omega(t,m)=A(t)\cap G(t)\cap G(v)\cap B(t,m)\cap K(t)\cap L(t)$  and
\begin{eqnarray*}
&A(t)=\left\{\max_{i\leq i_1} (K_{i+1}-K_i) \leq (\log (t))^2\right\},\\
&G(u)=\left\{\sup_{-u \leq r<s \leq  u } W_\kappa(s)-W_\kappa(r)\leq
\frac{1}{\kappa}( \log u  +3 \log\log  u)\right\},\\
&B(t,m)=\bigcap_{j=1}^{m-1}\left\{\sharp\{i\in N(- v, v):D_i
\geq \frac{1}{\kappa}\log { v^{k/m}} + 4 \log\log (v) \} \leq v^{1-\frac{k}{m}}\right\},\\
&K(t) =\left\{\sup_{\substack{-t<t_1<t_2<t\\ |t_2-t_1|<1}}|W_\kappa(t_2)-W_\kappa(t_1)|\leq (\log t)^{1/2}\log\log t\right\},\\
&L(t)=\left\{\sup_{0<r<s<v} W_{\kappa}(s)-W_{\kappa}(r) >\frac{1-\epsilon}{\kappa}\log v\right\}.
\end{eqnarray*}
Furthermore, whenever $u\rightarrow\infty$, the event $G(u)$ is fullfilled for $u$ large enough.
\end{lemma}
We now turn to some quenched estimates: let $[a,c]$ be an interval of $\mathbb{R}$. We call
\begin{equation}\label{defd}D_{+}=\sup_{x\in [a,c]}\left(\max_{y\in[x,c]}W_{\kappa}(y)-\min_{y\in[a,x)} W_{\kappa}(y)\right),\end{equation}
\begin{equation}D_{-}=\sup_{x\in [a,c]}\left(\max_{y\in[a,x]}W_{\kappa}(y)-\min_{y\in(x,c])} W_{\kappa}(y)\right),\end{equation}
and $$D=D_{-}\wedge D_{+}.$$
We also introduce $M:=\sup_{x\in [a,c]}W_{\kappa}(x)-\min_{x\in [a,c]}W_{\kappa}(x)$
We have
\begin{lemma}[quenched estimates]\label{qe}
Let $a,c,$ and $D$ be as above, and $W\in\Omega$, then for some constant $C$, and $u>1$
\begin{equation} \max_{x\in[a,c]}P_W^{x}\left[H(a) \wedge H(c)>Cu (M\vee 1)(1\vee(c-a)^4)) e^{D}\right]<e^{-u}.\label{tpssortie}\end{equation}
We also have a bound on the number of backtracks. For
$f\rightarrow \infty,$ $f=O(t)$
\begin{equation}P_W\left[\mathcal{B}\geq f \right]\leq C_3e^{-f}.\label{backtracks}\end{equation}
Finally, if $W\in \Omega$, for some constant $\gamma,$ for every $1\leq i \leq i_1$, and for $t$ large enough,
\begin{equation}\label{lpun}P_W^{K_i}\left[H(K_{i+1})>u\gamma (\log t)^{20} e^{D_{i-1}\vee
    D_{i}}|H(K_{i+1})<H(K_{i-1})\right]\leq e^{-u},\end{equation}
\begin{equation}\label{lpdeux}P_W^{K_i}\left[H(K_{i-1})>u\gamma (\log t)^{20} e^{D_{i-1}\vee
    D_{i}}|H(K_{i-1})<H(K_{i+1})\right]\leq e^{-u},\end{equation}
\begin{equation}\label{lptrois}P_W^{0}\left[H(K_{i_0})\wedge H(K_{i_0+1})>u\gamma (\log t)^{20} e^{D_{i_0-1}\vee
    D_{i_0}}\right]\leq e^{-u}.\end{equation}
\end{lemma}
Thanks to these lemmas, we are able to finish the proof of Theorem \ref{pqsl}.
\subsubsection{Upper bound.}
We recall $v=t^\nu$. Suppose $\Omega(t,m)$ is fulfilled, by the previous decomposition,
\begin{multline*}
P_{W}\left(H(v)>t\right)\leq P_{W}\left(H_{init}>\frac{t}{5} \right)+P_{W}\left(H_{dir}>\frac{t}{5} \right)\\+P_{W}\left(H_{back}>\frac{t}{5} \right)+P_{W}\left(H_{left}>\frac{t}{5} \right)+
P_{W}\left(H_{right}>\frac{t}{5} \right).
\end{multline*}
We begin with $H_{init}$. We recall from (\ref{hinit}) that $H_{init}$ is the time the diffusion takes to get to $K_{i_0+1}.$
Using the precedent estimates, on $G(v)$, we have, for $t$ large enough
$$D_{i_{0}}\vee D_{i_{0}+1}<\frac{1}{\kappa}(\log{v}+3\log\log v).$$
Thus, for every $\epsilon>0$,
\begin{multline*}P_{W}\left(H_{init}>\frac{t}{5}
  \right)\\ \leq P_{W}^0\left(H(K_{i_0+1})>\frac{t^{1-\nu/\kappa}}{5}e^{D_{i_{0}}\vee D_{i_{0}+1}}\cap H(K_{i_0+1})<H(K_{i_0})\right)\\+P_{W}^0\left(H(K_{i_0})>\frac{t^{1-\nu/\kappa}}{10}e^{D_{i_{0}}\vee D_{i_{0}+1}}\cap H(K_{i_0})<H(K_{i_0+1})\right)\\
  +P_W^{K_{i_0}}\left[H(K_{{i_0}+1})>\frac{t^{1-\nu/\kappa}}{10}e^{D_{i_{0}}\vee D_{i_{0}+1}}|H(K_{i_0+1})<H(K_{i_0-1})\right]
  \leq
 3 e^{-t^{1-\nu/\kappa-\epsilon}}.\\
\end{multline*}
Similarly, we have
\begin{multline*}P_{W}\left(H_{right}>\frac{t}{5}
  \right)=P_{W}^{K_{i_1}}\left(H(v)>\frac{t}{5}|H(v)<H(K_{i_1-1})\right)\leq
  e^{-t^{1-\nu/\kappa-\epsilon}}.\\
\end{multline*}
It is also a direct consequence of lemma \ref{qe} that, on $A(t),$ $i_0>\frac{t}{2(\log t)^2}$, whence, recalling the definition of $\mathcal{B}$ in (\ref{defb}),
\begin{equation*}
P_{W}\left(H_{left}>\frac{t}{5} \right)\leq
P_{W}\left(\mathcal{B}\geq\frac{t}{4\log^2t}\right)\leq\exp{\left(-\frac{t}{4\log^2 t}\right)}.\\
\end{equation*}

To deal with $H_{dir}$, note that
$$H_{dir}=\sum_{i=i_{0}+1}^{i_{1}-1} \tau_{+}^{(0)}(i),$$
where 
$\tau_+^{(0)}(i)$ is the first crossing of the interval $[K_i,K_{i+1}]$.
The $\tau_+^{(0)}(i)$ are independent random variables, and
$\tau_+^{(0)}(i)$ follows the same law as $H(K_{i+1})$ under
$P_W^{K_i}[\cdot|H(K_{i+1})<H(K_{i-1})]$.
 
On the other hand, if $H_{dir}>t/5,$ then the process spends an amount of time greater than $t/20m$ in the valleys of depth in
$$\left[\frac{ k}{\kappa m} \log{v} +4 \log\log {v},\frac{(k+1)}{\kappa m} \log{v} +4 \log\log {v}\right].$$
On $\Omega(t,m),$ the number of such valleys is at most $v^{1-\frac{k}{m}},$ we call $\sigma(k)$ the time spent in those valleys. By lemma \ref{qe}, and the precedent remarks, for some constant $C$,
$$\frac{\sigma(k)}{C(\log t)^{11} v^{(k+1)/\kappa m}} \vartriangleleft 2v^{(1-k/m)}+ \Gamma\left(2\lceil v^{(1-k/m)}\rceil,1\right),$$
where we note $A\vartriangleleft B$ for `` A is stochastically dominated by $B$'', and $\Gamma(k,\beta)$ is the Gamma distribution of parameter $(k,\beta)$.

For $m$ large enough, one can check easily that $\nu(1-k/m)<1-\nu(k+1)/m$ for all $k\leq m$,  whence, for $t$ large enough,

\begin{multline*}P_{W}\left[\sigma(k)\geq \frac{t}{20m}\right]\leq P\left[\Gamma\left(2v^{(1-k/m)},1\right)>\frac{t^{1-\nu(k+1)/\kappa m}}{(\log t)^{12}} \right] \\
\leq 4^{t^{\nu(1-k/m)}} \exp{\left(-\frac{t^{1-\nu(k+1)/\kappa m}}{(\log t)^{12}} \right)}\\
\leq \exp \left(-2t^{1-\nu(k+2)/(\kappa m)} +\log(4) t^{\nu(1-k/m)}\right).
\end{multline*}

Therefore, as $t\rightarrow \infty$,
$$P_{W}[H_{dir}>t/5]\leq m \exp \left(-t^{1-\nu(k+2)/(\kappa
  m)}\right) \leq m\exp\left(-t^{1-(1+\frac{2}{m})\frac{\nu}{\kappa
  }}\right).$$

We now deal with $H_{back}$.
\begin{multline*}
P_{W}\left(H_{back}>\frac{t}{5}
\right)\\ \leq\sum_{k=0}^{m-1}P_{W}\left(H_{back}>\frac{t}{5},\mathcal{B}\in[t^{k/m},t^{(k+1)/m}]\right) +P_W[\mathcal{B}>t].
\end{multline*}
By lemma \ref{qe}, 
$P_W[\mathcal{B}>t]<e^{-t},$ and 
\begin{equation}\label{inegalite}P_{W}\left(H_{back}>\frac{t}{5},\mathcal{B}\in[t^{k/m},t^{(k+1)/m}]\right)\leq
C\exp{\left(-t^{k/m}\right)}.\end{equation}

 On the other hand, 
$$H_{back}=\sum_{i=1}^{i_1-2} \sum_{j=1}^{\xi(i)}
 \tau_+^{(j)}(i)+\tau_-^{(j)}(i),$$
where
\begin{itemize}
\item $\tau_-^{(j)}(i)$ is the $j-th$ crossing of the interval  $[K_{i+1},K_i]$. 
\item $\tau_+^{(j)}(i)$ is the first crossing of the interval $[K_i,K_{i+1}]$ after the $j-th$ crossing of the interval  $[K_{i+1},K_i]$.  
\end{itemize}

The $\tau_{+,-}^{(j)}(i)$ are independent variables, and 
$\tau_+^{(j)}(i)$ follows the same law as $H(K_{i+1})$ under
$P_W^{K_i}[\cdot|H(K_i+1)<H(K_{i-1})],$ and $\tau_-^{(j)}(i)$ follows the same law as $H(K_{i})$ under
$P_W^{K_{i+1}}[\cdot|H(K_i)<H(K_{i+2})],$ (with the convention that
$K_{i_1+1}=v$). Therefore, thanks to lemma \ref{qe}, 

$$\frac{\tau_{+,-}^{(j)}(i)}{C e^H(\log t)^{10}}\vartriangleleft 1+{\bf e}$$
for some constant $C$ and $$H=\max_{i\in N\left(-\frac{t^{(k+1)}\log^2 t }{m},v\right)}D_i.$$
 Then, for
$W_\kappa \in \Omega(n,m)\cap G\left(\frac{t^{(k+1)}\log^2 t }{m}\right)$, on the event $\{\mathcal{B}\in[t^{k/m},t^{(k+1)/m}]\},$
$$\frac{H_{back}}{C (t^{(k+1)/m\kappa}\vee v^{1/\kappa})(\log t)^{10}}\vartriangleleft
2t^{(k+1)/m}+\Gamma(2t^{(k+1)/m},1).$$
Therefore, when $1-\frac{1}{\kappa}\left(\nu\vee \frac{k+1}{m}\right)\geq \frac{k+1}{m},$

$$P_{W}\left(H_{back}>\frac{t}{5},\mathcal{B}\in[t^{k/m},t^{(k+1)/m}]\right)\leq
C\exp\left(-C' t^{{1-\frac{1}{\kappa}\left(\nu\vee \frac{k+1}{m}\right)}} \right).$$
Putting this together with (\ref{inegalite}), we obtain

\begin{multline*}P_{W}\left(H_{back}>\frac{t}{5},\mathcal{B}\in[t^{k/m},t^{(k+1)/m}]\right) \\\leq
C\exp\left(-C' t^{{\left(1-\frac{1}{\kappa}\left(\nu\vee \frac{k+1}{m}\right)\right)} \vee \frac{k}{m}-\frac{1}{m}}\right).\end{multline*}

Putting together all the estimates, we get
\begin{multline*}\liminf_{t\rightarrow \infty} \frac {\log(-\log P_{W}
  [H(t^{\nu})>t])}{\log t}\\
\geq \min_{k\in [-1,m+1]} \left(\frac{k}{m}
  \vee \left( 1-\frac{1}{\kappa}\left(\nu\vee \frac{k+1}{m}\right)\right)-\frac{1}{m}\right) \wedge\left({1-(1+\frac{2}{m})\frac{\nu}{\kappa
  }}\right)
\\\geq \left(1-\frac{\nu}{\kappa}\right)\wedge \frac{\kappa}{\kappa
  +1}-\frac{3}{(1\wedge \kappa)m},\; P-a.s..
\end{multline*}
By taking the limit as $m$ goes to infinity, we get the upper bound for $P_{W}
  [H(t^{\nu})>t],$
namely
$$\liminf_{t\rightarrow \infty} \frac {\log(-\log P_{W}
  [H(t^{\nu})>t])}{\log t}\geq \left(1-\frac{\nu}{\kappa}\right)\wedge \frac{\kappa}{\kappa
  +1}.$$

We now turn to the proof of the lower bound.
\subsubsection{Lower bound.}\label{lb}

We suppose that $L(t)$ is fullfilled, therefore there is one valley of depth greater than $\frac{1-\epsilon}{\kappa}\log v$ before $v$. Let $b$ be the bottom of this valley, and $c$ such that
$b<c$ and $$W_{\kappa}(c)-W_{\kappa}(b)=\frac{1-\epsilon}{\kappa}\log v.$$
It is easy to see that 
$H(v)\geq H(c)-H(b)$, whence
$$P_{W} [H(t^\nu)>t]\geq P_{W}^{b} [H(c)>t].$$
We can suppose, without loss of generality, that $b=0$.
By the time change representation from the preliminary statements, under $P_W,$
$H(c)=T_\kappa(\sigma(A_\kappa(c))),$ where $\sigma(x)$ is the first
hitting time of $x$ by a brownian motion $B$. Therefore 
\begin{multline*}H(c)=\int_{0}^{\sigma(A_{\kappa}(c))} e^{-2W_{\kappa}(A_{\kappa}^{-1}(B_{s}))}
  d s.\\
=\int_{-\infty}^{A_{\kappa}(c)} \exp{(-2W_\kappa(A_{\kappa}^{-1}(x)))}L^{x}_{\sigma(A_\kappa(c))}d x\\
=\int_{-\infty}^c \exp{(-W_{\kappa}(u))} L^{A_{\kappa}(u)}_{\sigma(A_\kappa(c))}d u.
\end{multline*}
The last equality coming from a change of variable in the integral.
By a scaling argument, we get
$$H(c)\egall \int_{-\infty}^c \exp{(-W_{\kappa}(u))} A_{\kappa}(c) L^{A_{\kappa}(u)/A_{\kappa}(c)}_{\sigma(1)}d u.$$
We suppose $W_{\kappa}\in K(t),$ so
$$A_{\kappa}(c)\geq e^{W_{\kappa}(c)-(\log t)^{2/3}}>t^{(1-2\epsilon)\frac{\nu}{\kappa}},$$
and $A_\kappa(-1)>-e^{-(\log t)^{2/3}}.$
Hence 
$$H(c)\vartriangleright t^{(1-3\epsilon)\frac{\nu}{\kappa}}\inf_{x\in [A_{\kappa}(-1)/A_\kappa(c),0]}L^{x}_{\sigma(1)}.
$$
For $t$ large enough, $A_{\kappa}(-1)/A_\kappa(c)>-1/2$.
Therefore by the first Ray-Knight theorem (Statement \ref{rn1})
\begin{multline*}P_W^b[H(c)>t]\geq P\left[\inf_{x\in [ -1/2,0]}L^{x}_{\sigma(1)}>t^{1-\frac{\nu}{\kappa}+\varepsilon}\right]
\\
\geq P\left[Z'_1>2t^{1-\frac{\nu}{\kappa}+\varepsilon}\right]P\left[\sup_{u\in{[0,1/2]}}|Z_u|<t^{1-\frac{\nu}{\kappa}+\varepsilon}\right],
\end{multline*}
where $Z_t$ is a squared Bessel process of dimension $0$ started at $0$ and $Z'_t$ is a squared Bessel process of dimension $2$ started at $0$.
The last probability is greater than $1/2$ for $t$ large enough, and the first one is explicitly known (see 
for example \cite{Borodin:2002}).
We obtain that, for all $\varepsilon>0$, 
$$P_{W} [H(v)>t]\geq \exp\left(-\frac{1}{2}t^{1-\frac{\nu}{\kappa}+\varepsilon}\right).$$

To obtain the other lower bound, note that, similarly to lemma \ref{enve}, almost surely, there is a valley of depth at least 
$\frac{1-\epsilon}{\kappa+1}\log{t}$ in $[-t^{\kappa/(\kappa+1)},0]$, let $b'$ be the bottom of such valley, and $c'>b'$ such that $$W_{\kappa}(c')-W_{\kappa}(b')\geq\frac{1-\epsilon}{\kappa+1}\log{t}.$$
We have 
$$P_W[H(v)>t]\geq P_{W}[H(b)<H(t^{\nu})]P_{W}^b[H(c)>t].$$
Recalling the time change representation,
$$P_{W}[H(b)<H(t^{\nu})]=\frac{A_{\kappa}(t^\nu)}{A_{\kappa}(t^{\nu})-A_{\kappa}(b)}.$$
when $W_{\kappa}\in K(t)$, we can easily show that for every $\epsilon>0$, as $n$ goes to infinity, 
$$P_{W}[H(b)<H(t^{\nu})]\geq \exp{-t^{\frac{\kappa}{\kappa+1}+\epsilon}}.$$

By the same computations as for the first bound, we get
$$P_{W}^b[H(c)>t]\geq \exp{-t^{\frac{\kappa}{\kappa+1}+\epsilon}}.$$
Putting together both inequalities, we get  
$$\liminf_{t\rightarrow \infty} \frac {\log(-\log P_{W}
  [H(t^{\nu})>t])}{\log t}\leq \left(1-\frac{\nu}{\kappa}\right)\wedge \frac{\kappa}{\kappa
  +1},$$
which finishes the proof of Theorem \ref{pqsl}.
 
\subsection{Quenched slowdown for the diffusion.}
In this section we finish the proof of Theorem \ref{pqsl}.
The lower bound is trivial, since 
$$P_{W}[X_t<t^{\nu}]\geq P_{W}[H({t^{\nu}})>t].$$
To get the upper bound, let $m\in \mathbb{N}$, note that
\begin{multline}\label{decomposition}P_{W}[X_t<t^{\nu}]\leq P_{W}[H(t^{\nu})>t]
\\+\sum_{k=0}^{m-1} P_{W}\left[H\left(t^{\nu+\frac{k}{m}}\right)<t<H\left(t^{\nu+\frac{k+1}{m}}\right)\right]P_{W}^{t^{\nu+\frac{k}{m}}}[H(t^{\nu})<t]
\\+P_{W}[H(t^{\nu+1})<t]P_{W}^{t^{\nu+1}}[H(t^{\nu})<t].
 \\\end{multline}
Using the explicit distribution of the supremum before $t$ of a drifted brownian motion (see page 197 of \cite{Borodin:2002}) and the Borel-Cantelli lemma, we can easily see that for every $k\in \{1,m\},$ the event
$$U_m^k(n):= \left\{\sup_{(n+1)^{\nu}< s<t<n^{\nu+\frac{k}{m}}}W_{\kappa}(s)-W_{\kappa}(t) \geq\frac{\kappa}{4}n^{\nu+\frac{k}{m}}\right\}$$
is fullfilled for all $n$ large enough, therefore so does
$$U_n=\bigcup_{k=1}^m U_m^k(n).$$ 
Hence on $U_{\lceil t \rceil}$, there exist
$t^{\nu}<a<b<t^{\nu+\frac{k}{m}}$ such that
$$W_{\kappa}(a)-W_{\kappa}(b)\geq \frac{\kappa}{4}t^{\nu+\frac{k}{m}}.$$
By the same computations as in part \ref{lb}, we get that, on $U(\lceil t \rceil)$,
\begin{multline*}
 P_{W}^{b}[H(a)<t]
\leq P_W\left[e^{\frac{\kappa}{8}t^{\nu+\frac{k}{m}}}\inf_{x\in [0,e^{-\frac{\kappa}{8}t^{\nu+\frac{k}{m}}}]}L^{x}_{\sigma(1)}<t\right]\\
\leq
P_W\left[\inf_{u\in [1,1-e^{-\frac{\kappa}{8}t^{\nu+\frac{k}{m}}}]}Z_u<t e^{-\frac{\kappa}{8}t^{\nu+\frac{k}{m}}}\right],
\end{multline*}
where $Z_u$ is a squared Bessel process of dimension $2$ started at zero.
We have
\begin{multline*}P_W\left[\inf_{u\in [1,1-e^{-\frac{\kappa}{8}t^{\alpha}}]}Z_u<t e^{-\frac{\kappa}{8}t^{\alpha}}\right]\\\leq P_W\left[Z_1<2t e^{-\frac{\kappa}{8}t^{\alpha}}\right] +P_W\left[\sup_{u\in [1,1-e^{-\frac{\kappa}{8}t^{\alpha}}]}|Z_u-Z_1|\geq t e^{-\frac{\kappa}{8}t^{\alpha}}\right].\end{multline*}

Using statement \ref{Csorgo} with $u=t e^{-\frac{\kappa}{8}t^{\alpha}}$ and the fact that $\sqrt{Z_{1-t}-Z_1}$ is the Euclidean norm of a two dimensional Brownian motion, we get
$$P_W\left[\sup_{t\in [1,1-e^{-\frac{\kappa}{8}t^{\alpha}}]}|Z_t-Z_1|\geq t e^{-\frac{\kappa}{8}t^{\alpha}}\right]\leq2 \exp{-\frac{t}{10}}.$$
On the other hand, by the exact distribution of $Z_1$,
$$P(Z_1<x)=1-e^{-x/2}<x.$$
Therefore we get that for some constant $C$
$$ P_{W}^{t^{\nu+\frac{k}{m}}}[H(t^{\nu})<t]\leq P_{W}^{b}[H(a)<t]<e^{-C t^{\nu+\frac{k}{m}}}.$$

On the other hand, the bound for the hitting time implies that
$$P_{W}\left[H\left(t^{\nu+\frac{k}{m}}\right)<t<H\left(t^{\nu+\frac{k+1}{m}}\right)\right]\leq\exp{\left(-t^{\left(1-\left(\nu+\frac{k+1}{m}\right)/\kappa\right)\wedge\left(\frac{\kappa}{\kappa+1}\right)-\frac{1}{m}}\right)},$$
indeed the bound is trivial when $\nu+k/m>\kappa$.

The same arguments apply to the other terms of (\ref{decomposition}), whence
\begin{multline*}\liminf_{t\rightarrow\infty} \frac{\log(-\log P_{W}[X_t<t^{\nu}])}{\log t}
 \\\geq \min_{k\in [0,m]}\left[\left(\nu+\frac{k}{m}\right)\vee \left(\left(1-\frac{\nu+(k+1)/m}{\kappa}\right)\wedge\frac{\kappa}{\kappa+1} -\frac{1}{m}\right)\right].
\end{multline*}
Minimizing over $k$ and taking the limit as $m$ go to infinity, we get the desired upper bound.

\subsection{Proof of the lemmas.}
We begin with the estimates on the environment.
\subsubsection{Proof of lemma \ref{enve}.} 
Note that, as an easy consequence of statement \ref{Csorgo}, almost surely for $t$ large enough $i_1<2t.$ Therefore
\begin{equation}\label{at}A(t)\supset \tilde{A}(\lfloor t\rfloor):=\left\{\max_{i\leq2\lfloor t\rfloor+1}|K_{i+1}-K_i|\leq \log^2( \lfloor t\rfloor)\right\} .\end{equation}
Let us show that \begin{equation}{\cal P}[\tilde{A}(n)^c]=O(1/n^2)\label{an}.\end{equation} We have

\begin{equation}{\cal P}[\tilde{A}(n)^c]\leq \sum_{i=0}^{2n+1} {\cal P}[K_{i+1}-K_{i}\geq (\log{(n)})^2].\label{ansum}\end{equation}
By invariance of the environment, 
$${\cal P}[K_{1}-K_{0}\geq (\log{(n)})^2]={\cal P}[\tilde{K}_1\geq (\log{(n)})^2],$$
where
 $$\tilde{K}_{1}=\min\left\{t\geq 0: -\min_{s\in[0,n]}W_{\kappa}(s)\geq \frac{3}{\kappa} \log{n}, W_{\kappa}(t)>\sup_{s\geq t}W_{\kappa}(s)-1\right\}.$$
On the other hand, conditionally to $K_i$, the process $ W_{\kappa}(K_i+s)-W_{\kappa}(K_i)$ is a drifted Brownian motion conditionned to have its supremum lesser than $1$. Therefore
\begin{multline*}{\cal P}[K_{i+1}- K_{i}\geq (\log{n})^2]={\cal P}[\tilde{K}_{1}\geq(\log n)^2|\sup_{t\geq 0}W_{\kappa}\leq 1]\\
\leq\frac{{\cal P}[\tilde{K}_{1}\geq(\log n)^2]}{{\cal P}[\sup_{t\geq 0}W_{\kappa}\leq 1]}.
\end{multline*}
For $\kappa>0$, ${\cal P}[\sup_{t\geq 0}W_{\kappa}\leq 1]$ is a positive constant.
It remains to bound ${\cal P}[\tilde{K}_{1}\geq(\log n)^2]$, note that if $$W_\kappa\left((\log n)^2\right)<-\frac{6}{\kappa} \log{n},$$ and $$\sup_{t\geq (\log{n})^2} W_{\kappa}(t)-W_{\kappa}\left((\log n)^2\right)<\frac{3}{\kappa} \log{n},$$ then there exists one point $x^*$ before $(\log n)^2$ such that $\inf_{t\in[0,x^*]} W_{\kappa}(t)<-\frac{3}{\kappa}\log{n}$ and $W_{\kappa}(x^*)\geq \sup_{s\geq x^*}W_{\kappa}(s)-1$ (see figure 3), therefore $\tilde{K}_1<(\log n)^2$. Taking the complementary events, we get
\begin{multline*}{\cal P}[\tilde{K}_{1}\geq(\log n)^2]\\ \leq {\cal P}\left[W_{\kappa}\left((\log{n})^2\right)>-\frac{6}{\kappa} \log{n} \text{ or } \sup_{t\geq (\log{n})^2} W_{\kappa}(t)-W_{\kappa}\left((\log n)^2\right)>\frac{3}{\kappa} \log{n}\right].
\end{multline*}
\begin{figure}[h]
\centering
\scalebox{1} 
{
\begin{pspicture}(1,-2.06)(7.48,2.1)
\rput(2.2,0.94){\psaxes[linewidth=0.025999999,labels=none,ticks=none,ticksize=0.10583333cm]{->}(0,0)(0,-3)(5,1)}
\psline[linewidth=0.018cm,linestyle=dashed,dash=0.16cm 0.16cm](2.2,-0.22)(7.32,-0.22)
\psline[linewidth=0.018cm,linestyle=dashed,dash=0.16cm 0.16cm](2.2,-1.4)(7.32,-1.4)
\psline[linewidth=0.018cm](5.54,1.06)(5.54,0.82)
\rput(7.2,0.68){\footnotesize $x$}

\rput(2.89,1.92){\footnotesize $W_{\kappa}(x)$}

\rput(5.6,1.32){\footnotesize $(\log n)^2$}

\rput(1.47,-0.2){\footnotesize $-\frac{3}{\kappa}\log n$}

\rput(1.47,-1.38){\footnotesize $-\frac{6}{\kappa}\log n$}
\pscustom[linewidth=0.018]
{
\newpath
\moveto(2.2,0.92)
\lineto(2.24,0.89)
\curveto(2.26,0.875)(2.32,0.845)(2.36,0.83)
\curveto(2.4,0.815)(2.445,0.775)(2.45,0.75)
\curveto(2.455,0.725)(2.47,0.675)(2.48,0.65)
\curveto(2.49,0.625)(2.52,0.575)(2.54,0.55)
\curveto(2.56,0.525)(2.615,0.48)(2.65,0.46)
\curveto(2.685,0.44)(2.75,0.41)(2.78,0.4)
\curveto(2.81,0.39)(2.855,0.35)(2.87,0.32)
\curveto(2.885,0.29)(2.91,0.24)(2.92,0.22)
\curveto(2.93,0.2)(2.955,0.155)(2.97,0.13)
\curveto(2.985,0.105)(3.06,0.045)(3.12,0.01)
\curveto(3.18,-0.025)(3.3,-0.065)(3.36,-0.07)
\curveto(3.42,-0.075)(3.515,-0.095)(3.55,-0.11)
\curveto(3.585,-0.125)(3.62,-0.165)(3.62,-0.19)
\curveto(3.62,-0.215)(3.64,-0.265)(3.66,-0.29)
\curveto(3.68,-0.315)(3.755,-0.35)(3.81,-0.36)
\curveto(3.865,-0.37)(3.94,-0.37)(3.96,-0.36)
\curveto(3.98,-0.35)(4.01,-0.32)(4.02,-0.3)
\curveto(4.03,-0.28)(4.06,-0.24)(4.08,-0.22)
\curveto(4.1,-0.2)(4.145,-0.17)(4.17,-0.16)
\curveto(4.195,-0.15)(4.24,-0.14)(4.26,-0.14)
\curveto(4.28,-0.14)(4.325,-0.18)(4.35,-0.22)
\curveto(4.375,-0.26)(4.415,-0.32)(4.43,-0.34)
\curveto(4.445,-0.36)(4.49,-0.42)(4.52,-0.46)
\curveto(4.55,-0.5)(4.62,-0.565)(4.66,-0.59)
\curveto(4.7,-0.615)(4.76,-0.66)(4.78,-0.68)
\curveto(4.8,-0.7)(4.84,-0.775)(4.86,-0.83)
\curveto(4.88,-0.885)(4.905,-0.98)(4.91,-1.02)
\curveto(4.915,-1.06)(4.925,-1.155)(4.93,-1.21)
\curveto(4.935,-1.265)(4.96,-1.345)(4.98,-1.37)
\curveto(5.0,-1.395)(5.055,-1.445)(5.09,-1.47)
\curveto(5.125,-1.495)(5.19,-1.53)(5.22,-1.54)
\curveto(5.25,-1.55)(5.32,-1.56)(5.36,-1.56)
\curveto(5.4,-1.56)(5.45,-1.575)(5.46,-1.59)
\curveto(5.47,-1.605)(5.495,-1.64)(5.51,-1.66)
\curveto(5.525,-1.68)(5.575,-1.705)(5.61,-1.71)
\curveto(5.645,-1.715)(5.715,-1.71)(5.75,-1.7)
\curveto(5.785,-1.69)(5.85,-1.655)(5.88,-1.63)
\curveto(5.91,-1.605)(5.96,-1.545)(5.98,-1.51)
\curveto(6.0,-1.475)(6.045,-1.42)(6.07,-1.4)
\curveto(6.095,-1.38)(6.145,-1.35)(6.17,-1.34)
\curveto(6.195,-1.33)(6.245,-1.32)(6.27,-1.32)
\curveto(6.295,-1.32)(6.365,-1.31)(6.41,-1.3)
\curveto(6.455,-1.29)(6.525,-1.255)(6.55,-1.23)
\curveto(6.575,-1.205)(6.625,-1.145)(6.65,-1.11)
\curveto(6.675,-1.075)(6.725,-1.005)(6.75,-0.97)
\curveto(6.775,-0.935)(6.825,-0.905)(6.85,-0.91)
\curveto(6.875,-0.915)(6.93,-0.955)(6.96,-0.99)
\curveto(6.99,-1.025)(7.045,-1.06)(7.07,-1.06)
\curveto(7.095,-1.06)(7.13,-1.06)(7.16,-1.06)
}
\psline[linewidth=0.018cm,linestyle=dashed,dash=0.16cm 0.16cm](4.24,1.04)(4.24,-0.14)
\rput(4.29,1.32){\footnotesize $x^*$}
\end{pspicture}
} 

 \caption{$\tilde{K}_1$}
\end{figure}
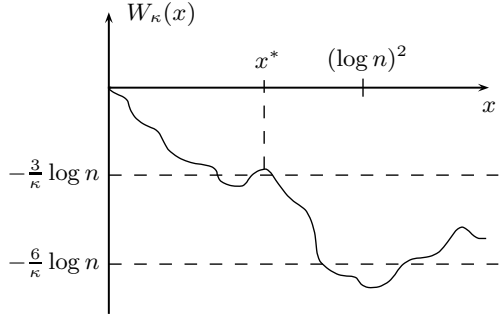

By standard gaussian estimates,
$${\cal P}\left[W_{\kappa}\left((\log{n})^2\right)>-\frac{6}{\kappa} \log{n}\right]=O(n^{-3})$$  and 
$${\cal P}\left[\sup_{t\geq (\log{n})^2} W_{\kappa}(t)-W_{\kappa}\left((\log{n})^2\right)>\frac{3}{\kappa} \log{n}\right]={\cal P}\left[\sup_{t\geq 0} W_{\kappa}(t)>\frac{3}{\kappa} \log{n}\right].$$
By formula 1.1.4(1) from page 197 of \cite{Borodin:2002},
the last probability is equal to $n^{-3}$. Therefore recalling equation (\ref{ansum}), this finishes the proof of (\ref{an}).
Therefore, using the Borel-Cantelli lemma and (\ref{at}), $A(t)$ is fullfilled for every $t$ large enough.

We now turn to $G$. We consider the process 
\begin{equation}\label{defut}U_{t}:= \sup_{-\infty\leq s \leq t} W_{\kappa}(t)-W_{\kappa}(s).\end{equation}
Note that for $n=\lfloor t\rfloor$,$$\left\{\sup_{-(n+1)\leq t\leq n+1}U_{t}\leq \frac{1}{\kappa}(\log
n + 3\log\log n)\right\}\subset G(t).$$

The process $U_{t}$ is called a Reflected Brownian Motion with drift. This kind of process appears naturally in some queueing system
models.
 It is a positive and stationnary diffusion, with stationnary
law the exponential law of parameter $\kappa$. It is also reversible in time, therefore we can reduce to
proving that, as $n$ goes to infinity, the event \begin{equation}\label{suput}\left\{\sup_{0\leq t\leq n+1}U_{t}\leq \frac{1}{\kappa}(\log
n + 3\log\log n)\right\}\end{equation} is fullfilled.

In \cite{Salminen:2001} it is shown that the length of the excursions
away from zero (or busy periods) of $U_{t}$ follows a gamma distribution $\Gamma\left(\frac{1}{2},\frac{\kappa^2}{8}\right),$ and that the supremum $m_{0}$ over one excursion of $U_t$ has an explicit law, given by
\begin{equation}{\cal P}(m_{0}>y)=\frac{2e^{-\kappa y}}{(1-e^{-\kappa y})^2}(\kappa y -(1-e^{-\kappa y})).\label{mdet}\end{equation}
Let $C$ be some large constant.
We call $F(n)$ the event that $U_{t}$ makes more than $C n$ excursions between time $0$ and time $n+1$.
We have 
$${\cal P}(F(n))\leq {\cal P}\left(\Gamma\left(\frac{C n}{2},\frac{\kappa^2}{8}\right)<n+1\right)=\frac{\gamma(C n/2,\frac{(n+1)\kappa^2}{8})}{\Gamma(C n/2)},$$
where $\gamma(\cdot,\cdot)$ is the incomplete gamma function.
By Stirling's formula, $${\cal P}(F(n))=O(((n+1)\kappa^2/8)^{C n/2}(C n/2e)^{-C n/2-1/2})=o(n^{-4})$$ for $C$ large enough.
Therefore by the Borel-Cantelli lemma, almost surely there exists $n_0$ such that $F(n)$ is fullfilled for all $n\geq n_0$.

On the other hand, we call $\tilde{G}(k)$ the event that the maximum during
the $k-th$ excursion is lower than $1/\kappa (\log k + 3 \log \log k ).$ Recalling (\ref{mdet}), for $k\geq 10,$
\begin{equation*}{\cal P}\left(\tilde{G}(k)^c\right)={\cal P}\left(m_{0}>\frac{1}{\kappa} \left(\log k +3 \log \log k \right)\right)\leq \frac{8}{k(\log k)^2}.\end{equation*}
By the Borel-Cantelli lemma, we get that there exists $k_0$ such that $\tilde{G}(k)$ is fullfilled for all
$k\geq k_0$. Take $n>n_0\vee k_0$, and such that 
$$\frac{1}{\kappa}(\log n + 3\log\log n)$$ is greater than the supremum over the $k_0$ first excursions of $U_t$.
Then on $F(n)\cap \bigcap_{k=k_0}^{n}\tilde{G}(k)$ the event in (\ref{suput}) is fullfilled. This implies the result for $G(t)$.
\\

Let us turn to $B(t,m)$.
Let $n=\lfloor v \rfloor$.
We call, for $0<a<1$ 
$$\tilde{B}(n,a)=\left\{\sharp\left\{i\in
  N[-(n+1),n+1]:D_i\geq\frac{a}{\kappa}\log n + \frac{4}{\kappa} \log\log n\right\}<n^{1-a}\right\}.$$
Recalling the definitions of the $K_i$ and $U_t$, we note that the event that
two different $K_i$ belong to the same excursion of $U_t$ implies that
 the maximum during this excursion is at least $3/\kappa \log n$,
  therefore, by the same argument as before, when $n$ is large enough, this does not happen. We can also suppose that $U_t$ makes less than $C n$ excursions between
  time $-(n+1)$ and $n+1$.
Thus, on these events, 
$$\sharp\left\{i\in
  N[-(n+1),(n+1)]:H_i\geq\frac{a}{\kappa}\log n + 4
  \log\log n\right\}$$ is stochastically dominated by a
  $Binomial(2n+1,p)$, where 
$$p=P\left[m_t \geq \frac{a}{\kappa}\log
    n+4\log\log n \right]<2\frac{n^{-a}}{\log n ^{2}}.$$
Whence, using Chebyshev's  exponential inequality, 
\begin{multline*}P[\tilde{B}(n,a)^c]\leq
\exp{\left(-n^{1-a}\right)}\exp{\left((2n+1)\log(1+p(e-1))\right)}\\\leq\exp\left(4n p-n^{1-a}\right).\end{multline*}
The estimate on $p$, together with the Borel-Cantelli lemma, implies that, almost surely for $n$ large enough, 
$$\bigcap_1^{m-1}\tilde{B}(n,k/m)\subset B(t,m)$$ is fullfilled.
\\

We finally prove that $L(t)$ is fullfilled for $t$ large enough.
Recalling the notations concerning $U_t$ from (\ref{defut}), we call $f(n)$ the event that $U_t$ makes more that $\frac{n}{(\log n)^2}$ excursions before time $n$. Using the explicit distribution of the length of the excursions of $U_t$, we have
$${\cal P}\left(f(n)^c\right)\leq  {\cal P}\left(\Gamma\left(\frac{ n}{2{{(\log n)^2}}},\frac{\kappa^2}{8}\right)>n\right).$$
Recalling that a $\Gamma(k,\theta)$ distribution has expectation $k\theta$ and variance $k\theta^2$, by Bienaym\'e-Chebyshev's inequality, for $n$ large, $${\cal P}\left(f(n)^c\right)\leq\frac{10}{n(\log n)^2}.$$
Now the Borel-Cantelli lemma implies that $f(n)$ is fullfilled for all $n$ large enough.

Now suppose that $f(\lfloor v \rfloor)$ is fullfilled,
Note that $U_t$ and $\sup_{0<s<t} W_{\kappa}(t)-W_{\kappa}(s)$ are equal after the first $0$ of $U_t$. Call $\tilde{L}(t)$ the event that there exists one excursion of height at least $\frac{1-\epsilon}{\kappa}\log (n+1)$ between the second and the $\lfloor \frac{n}{(\log n)^2}\rfloor$-th excursion of $U_t$. It is easy to see that 
$$f(\lfloor v \rfloor)\cap\tilde{L}(t)\subset L(t).$$

On the other hand, by (\ref{mdet}),
\begin{multline*}{\cal P}\left(\tilde{L}(t)^c\right)\leq {\cal P}\left[m_t<\frac{1-\epsilon}{\kappa}\log(n+1)\right]^{\frac{n}{(\log{n})^2}}\\\leq\left(1-e^{-(1-\epsilon)\log{n+1}}\right)^{\frac{n}{(\log{n})^2}}.\end{multline*}
This is summable, therefore we can apply the Borel-Cantelli lemma to get the result on $\tilde{L}(t)$, then on $L(t)$.\\

The result on $K(t)$ is a direct consequence of statement \ref{Csorgo}.\\

We now turn to the quenched estimates.
\subsubsection{Proof of Lemma \ref{qe}. }
We begin with the proof of (\ref{tpssortie}). Without loss of generality we can suppose $x=0$ and $D=D_{+}.$ 
We suppose $|c-a|\geq1$, the proof being similar when $|c-a|\leq1.$

Recalling from the preliminary statements the time change representation of
$X_t$, 
 we get that, under $P_W,$
$H(v)=T_\kappa(\sigma(A_\kappa(v))),$ where 
$$A_\kappa(x)=\int_0^x e^{W\kappa(y)}d y,$$
$$T_\kappa(t)=\int_0^t e^{-2W_\kappa(A_\kappa^{-1}(B(s)))}d s,$$ and $\sigma(x)$ is the first
hitting time of $x$ by a Brownian motion $B.$ Therefore 
\begin{multline*}H(a)\wedge H(c)=\int_{0}^{\sigma(A_{\kappa}(a))\wedge
    \sigma(A_{\kappa}(c))} e^{-2W_{\kappa}(A_{\kappa}^{-1}(B_{s}))}
  d s\\
=\int_{A_{\kappa}(a)}^{A_{\kappa}(c)} \exp{(-2W_\kappa(A_\kappa^{-1}(x)))}L^{x}_{\sigma(A_\kappa(a))\wedge\sigma(A_\kappa(c))}d x.
\end{multline*}
We are going to use the second Ray-Knight Theorem (Statement \ref{rn}) : note that
$$L^x_{\sigma(A_\kappa(a))\wedge\sigma(A_\kappa(c))}\leq
L^x_{\sigma(A_\kappa(c))},$$ and that $L^x_{\sigma(A_\kappa(c))}$ is
stochastically dominated by the local time at $x$ before
$\sigma(A_\kappa(c))$ of a Brownian motion started at $a$.
Therefore
$$H(a)\wedge H(b)\vartriangleleft
\int_0^{A_\kappa(c)-A_{\kappa}(a)}{V}(s)X_s d s,$$
where ${V}(x)= \exp{(-2W_\kappa(A_\kappa^{-1}(A_\kappa(c)-x)))},$ and $X_{s}$ is a Bessel process of dimension 2, started at 0. We
call $\alpha:=A_\kappa(c)-A_{\kappa}(a)$, and $\lambda(V)$ the supremum of all $\lambda$ such that a
solution to 
$$y''(t)=-\lambda V(t)y(t),\;t\geq0\;y'(\alpha)=0,\;y(\alpha)=1$$ is positive
in $[0,\alpha]$.
$\lambda(V)$ is usually known as the spectral gap, or Poincar\'e's constant
associated to V.

By a standard change of variable in the previous differential
equation, and an application of statement \ref{bobkov}, we get
\begin{multline*}\frac{1}{\lambda(V)}\\ \leq 32 (A_{\kappa}(c)-A_{\kappa}(a))^2 \sup_{0<t<1} (1-t) \int_{0}^t e^{-2W_{\kappa}( A_{\kappa}^{-1}\left(A_{\kappa}(a)+s(A_{\kappa}(c)-A_{\kappa}(a))  \right))}d s\\
=32 (A_{\kappa}(c)-A_{\kappa}(a)) \sup_{0<t<1} (1-t)\int_{A_{\kappa}(a)}^
{A_{\kappa}(a)+t (A_{\kappa}(c)-A_{\kappa}(a))} e^{-2W_{\kappa}(A_{\kappa}^{-1}(u))} d u\\
=32 (A_{\kappa}(c)-A_{\kappa}(a)) \sup_{0<t<1} (1-t) \int_{a}^{d(t)} e^{-W_{\kappa}(v)}d v,
\end{multline*}
where $d(t)=A_\kappa^{-1}({A_{\kappa}(a)+t
  (A_{\kappa}(c)-A_{\kappa}(a))})$.
Easy computations show that 
$$(1-t)(A_{\kappa}(c)-A_{\kappa}(a))= \int_{d(t)}^{c}
e^{W_{\kappa}(v)}d v,$$ whence, recalling from (\ref{defd}) that $$D_{+}=\sup_{x\in [a,c]}\left(\max_{y\in[x,c]}W_{\kappa}(y)-\min_{y\in[a,x)} W_{\kappa}(y)\right),$$ we get
$$\frac{1}{\lambda(V)} \leq 32\sup_{a\leq x\leq c}\int_{a}^{x}
e^{-W_{\kappa}(v)}d v \int_{x}^{c}e^{W_{\kappa}(v)}d v\leq 32 (c-a)e^{D^+}.$$

From Lemma \ref{Laplace} we get that $E[\exp{\lambda(V)U}]$ is finite,
but we need an explicit bound. Toward this goal we are going to extend
the interval : let $c'$ be such that $(c'-a)=2(c-a)$ and let us extend $W_\kappa$ on
$[c,c']$ by a constant function (equal to $W_{\kappa}(c)$).
We call ${\tilde V}(x)= \exp{(-2W_\kappa(A_\kappa^{-1}(A_\kappa(c)-x)))},$ for $x\in[A_{\kappa}(c)-A_{\kappa}(c'),A_{\kappa}(c)-A_{\kappa}(a)]$ and $\lambda(\tilde{V})$ the supremum of all $\lambda$ such that a
solution to 
\begin{equation}\label{equadiff}y''(t)=-\lambda \tilde{V}(t)y(t),\;t\geq\;y'(\alpha)=0,\;y(\alpha)=1\end{equation} is positive
in $[A_{\kappa}(c)-A_{\kappa}(c'),\alpha]$.

By the same calculations as before we get 
   \begin{multline*}\frac{1}{\lambda(\tilde{V})} \leq 32\sup_{a\leq x\leq c'}\int_{a}^{x}
e^{-W_{\kappa}(v)}d v \int_{x}^{c'}e^{W_{\kappa}(v)}d v\leq 32 (c'-a)e^{D^+}\\=64 (c-a)e^{D^+}.\end{multline*}
For $\lambda<\lambda(\tilde{V})$, let $\phi$ be a solution to (\ref{equadiff}) on $[A_{\kappa}(c)-A_{\kappa}(c'),\alpha]$, then $\phi$ is a solution to (\ref{equadiff}) on $[0,\alpha]$, and by concavity, 
$$\phi(0)\geq \frac{A_{\kappa}(c')-A_{\kappa}(c)}{A_{\kappa}(c')-A_{\kappa}(a)}\geq \frac{e^{-M}}{2}.$$
Together with lemma \ref{Laplace}, we get $$E_{W}[\exp(\lambda H(a)\wedge H(c))]<2e^M.$$
This, together with Markov's inequality, finishes the proof of the first part of lemma \ref{qe}.
\\

In order to prove (\ref{backtracks}), note that, due to the time change 
representation, and for $W\in \Omega$,
\begin{multline}\label{ret} P_W^{K_i}\left[H({K_{i-1}})<H({K_{i+1}})\right]=\int_{K_i}^{K_{i+1}}e^{W_\kappa(x)}d x
\left(\int_{K_{i-1}}^{K_{i+1}}e^{W_\kappa(x)}d x\right)^{-1}\\
\leq \max_{i\leq
  i_1}(K_i-K_{i-1})\frac{e^{1+{W_\kappa(K_i)}}}{e^{W_{\kappa}(K_{i-1})-(\log{t})^{1/2}\log\log t}}
\leq t^{-3/2},
\end{multline}
using the fact, that, by definition of the $K_i$, on $K(t)\cap G(t)$, $$W_{\kappa}(K_{i-1})\geq\inf_{K_{i-1}\leq x
  \leq K_i} W_\kappa(x)+\frac{3}{\kappa} \log t\geq W_{\kappa}(K_i)+\frac{2}{\kappa}\log t -\frac{3}{\kappa}\log\log t
  .$$
Then we have to distinguish two cases : either the walk $Y_j$ gets to
  the level $v$ in more than $3n$ steps or in less than $3n$ steps.
In the first case there are at least $n$ steps back before $H(v)$, and
  in the second case the number of steps back is dominated by a
  $Binomial(3n,n^{-3/2})$. Thus
$$P_W\left[\mathcal{B}\geq f(t)\right]\leq\left(\begin{array}{c} 3n\\n
  \end{array}\right) \left( \frac{1}{n^{3/2}}\right)^n +
  P\left[Binomial(3n,n^{-3/2})\geq f(t)\right].$$
The result follows easily from Stirling's formula and Chebyshev's
exponential inequality.\\

We now turn to the proof of (\ref{lpun}),(\ref{lpdeux}) and (\ref{lptrois}). We start with (\ref{lpun}).
First note that 
\begin{multline*}P_W^{K_i}\left[H(K_{i+1})>u\gamma (\log t)^{20} e^{D_{i-1}\vee
    D_{i}}|H(K_{i+1})<H(K_{i-1})\right]\\\leq \frac{P_W^{K_i}\left[H(K_{i-1})\wedge H(K_{i+1})>u\gamma (\log t)^{20} e^{D_{i-1}\vee
    D_{i}})\right]}{P_W^{K_i}\left[H(K_{i+1})<H(K_{i-1})\right]},\end{multline*}
As a direct consequence of (\ref{ret}), we have, $P-a.s.,$ for $n$ large enough,
$$P_W^{K_i}\left[H({K_{i+1}})<H({K_{i-1}})\right]\geq \frac{1}{2}.$$

We are going to use (\ref{tpssortie}) in order to bound the numerator.
Note that, due to the definition of the $K_{i},$ 
$$\sup_{K_{i-1}<s<t<K_{i+1}}W_{\kappa}(s)-W_{\kappa}(t)\geq D_{i-1}\vee D_{i}. $$
On the other hand, on $A(t)\cap K(t)$, 
$$K_{i+1}-K_{i-1}\leq 2(\log t)^2,$$ and then
$$\sup_{x\in [K_{i-1},K_{i+1}]}W_{\kappa}(x)-\min_{x\in [K_{i-1},K_{i+1}]}W_{\kappa}(x)<(\log t )^3.$$
Therefore, the result follows easily by application of (\ref{tpssortie}).

We now turn to the proof of (\ref{lpdeux}).
As before,
 \begin{multline*}P_W^{K_i}\left[H(K_{i-1})>u\gamma (\log t)^{20} e^{D_{i-1}\vee
    D_{i}}|H(K_{i-1})<H(K_{i+1})\right]\\\leq \frac{P_W^{K_i}\left[H(K_{i-1})\wedge H(K_{i+1})>u\gamma (\log t)^{20} e^{D_{i-1}\vee
    D_{i}})\right]}{P_W^{K_i}\left[H(K_{i-1})<H(K_{i+1})\right]}.\end{multline*}
The numerator is the same as in the proof of (\ref{lpun}), so we only have to deal with the denominator.
We recall from (\ref{ret}) that $$P_W^{K_i}\left[H({K_{i-1}})<H({K_{i+1}})\right]=\int_{K_i}^{K_{i+1}}e^{W_\kappa(x)}d x
\left(\int_{K_{i-1}}^{K_{i+1}}e^{W_\kappa(x)}d x\right)^{-1}.$$ On $K(t)\cap G(t)$, we obtain easily
$$P_W^{K_i}\left[H({K_{i-1}})<H({K_{i+1}})\right]\geq \frac{e^{W_{\kappa}(K_i)-W_{\kappa}(K_{i-1})-\log t}}{\left(\log t\right)^3}.$$
Note that on $A(t)\cap K(t)$, $W_{\kappa}(K_{i-1})-W_{\kappa}(K_i)\leq (\log t)^3.$ (\ref{lpdeux}) follows then easily.

The proof of (\ref{lptrois}) is similar and omitted.
\section{Quenched speedup.}
In this part we show Theorem \ref{qs}. We first recall some facts.
\subsection{Preliminary statements.}
Our proof is mainly based on ``Kotani's formula'', expressed in
\cite{Kawazu:1997},
\begin{statement}[Kotani's lemma]
Let $\lambda>0$. Then for $t\geq 0$
$$E_W\left[e^{-\lambda H(t)}\right]= \exp{\left(-2\lambda \int_0^t
  U_\lambda(s) d s\right)}, \;P-a.s.,$$
where $U_\lambda(t)$ is the unique stationnary and positive solution
  of the equation
$$d U_\lambda(t)=U_\lambda(t) dW(t)+\left(1 +\frac{1-\kappa}{2}
  U_\lambda(t)- 2\lambda U_\lambda(t)^2\right)d t.$$
(Here $W(t)$ is the Brownian motion defined in the introduction). 
\end{statement}
We shall also use the following result from \cite{Hu:2004} (Lemma 2.4)
\begin{statement}\label{tl}
$$\lim\frac{1}{r} \sup_{|x|<u}\left(L^x_{\tau(r)} - r\right)=0,\; a.s.,$$
whenever $u\rightarrow \infty$ and $r\gg u\log\log u$.
\end{statement}

\subsection{Proof of Theorem \ref{pqs}.}
We use the same time change method as in the annealed case, in order
to get almost sure estimates for $U_\lambda$.
Let $$g(x)=\int_1^x \frac{e^{2/s + 4\lambda s}}{s^{1-\kappa}}d s.$$ One can easily check that $g$ is a
  scale function of $U_\lambda.$
  By the same arguments as in section \ref{tpp}, we get
  $$\int_{0}^t U_{\lambda}(s) d s= \int_{0}^{\mu(t)} g^{-1}(\gamma(u))^{1-2\kappa}\exp{\left(-\frac{4}{g^{-1}(\gamma(u))}-8\lambda g^{-1}(\gamma(u))\right)}d u,$$
where $\gamma(u)$ is a standard brownian motion,
$$\mu(t)=\int_{0}^t U_{\lambda}(s)^{2\kappa}\exp{\left(\frac{4}{U_{\lambda}(s)}+8\lambda U_{\lambda}(s)\right)}d s, $$
and $$\mu^{-1}(t)= \int_{0}^t
g^{-1}(\gamma(s))^{-2\kappa}\exp{\left(-\frac{4}{g^{-1}(\gamma(s))}-8\lambda
  g^{-1}(\gamma(s))\right)}d s.$$
We have the following lemma, whose proof is postponed
\begin{lemma}\label{equ}
Let $\nu\in \mathbb{R}$, and $$D_\nu(r)=\int_{0}^{\tau_r}
g^{-1}(\gamma(s))^{\nu}\exp{\left(-\frac{4}{g^{-1}(\gamma(s))}-8\lambda
  g^{-1}(\gamma((s))\right)}d s.$$
Then, whenever $\lambda\rightarrow 0$ and
$r\gg\log(1/\lambda)\log\log(1/\lambda),$
$$D_{1-2\kappa}(r)=r(1+o(1))\frac{\Gamma(1-
  \kappa)}{(4\lambda)^{1-\kappa}};$$ and for some positive constant
  $D$, 
$$D_{-2\kappa}(r)= D r(1+o(1)).$$
\end{lemma}
Let us use this lemma to finish the proof of Theorem \ref{pqs}.
We get easily that 
$\mu^{-1}(\tau_r)=D_{-2\kappa}(r).$
Whence, for some constant $D'$, $$\tau_{D'(1-o(1))t}\leq\mu(t) \leq \tau_{D'(1+o(1)) t}$$ almost
surely, as $\lambda\rightarrow 0$ and
$t\gg\log(1/\lambda)\log\log(1/\lambda).$
Therefore, under the same assumptions, for some constant $D''$, 
$$D''(1-o(1)){t}\frac{\Gamma(1-
  \kappa)}{(4\lambda)^{1-\kappa}}\leq \int_{0}^t U_{\lambda}(s) d s\leq  D''(1+o(1))t\frac{\Gamma(1-
  \kappa)}{(4\lambda)^{1-\kappa}}.$$
Thus, going back to Kotani's lemma, for $t>0$, and for some constant $C$, we get, as
  $\lambda\rightarrow 0$, $v\gg\log(1/\lambda)\log\log\log(1/\lambda),$
 \begin{equation} \label{hdev} \exp{\left(-C(1+o(1))\lambda^\kappa v \right)}\leq E_W\left[e^{-\lambda H(v)}\right]\leq
 \exp{\left(-{C}(1-o(1))\lambda^\kappa v \right)}, \;P-a.s..\end{equation}
 By application of Chebyshev's inequality, for $\lambda$ as before,
$$\log P_W\left[H(v)<\left(\frac{v}{u}\right)^{1/\kappa}\right] \leq \lambda \left(\frac{v}{u}\right)^{1/\kappa}-C(1-o(1)) v\lambda^\kappa.$$
 We call $\lambda(x)$ the value of lambda that minimizes $\lambda x-C v\lambda^\kappa.$ It is clear that $\lambda(x)$ is a decreasing function of $x$, such that
\begin{equation}\label{lambd}\lambda(x) x =C v\kappa \lambda(x)^{\kappa}.\end{equation}
Let $\lambda^*=\lambda\left(\left(\frac{v}{u}\right)^{1/\kappa}\right),$ we get easily the expression 
\begin{equation}\label{lambda}{\lambda^*}^\kappa = \left({C\kappa}\right)^{\frac{\kappa}{1-\kappa}} \frac{u^{\frac{1}{1-\kappa}}}{v}.\end{equation} 
One can easily check that $\lambda^*\rightarrow 0$, $v\gg\log(1/\lambda^*)\log\log\log(1/\lambda^*).$ Therefore we can apply the precedent estimate to get
 $$\limsup_{v\rightarrow \infty}\frac{\log P_W\left[H(v)<\left(\frac{v}{u}\right)^{1/\kappa}\right]}{u^{\frac{1}{1-\kappa}}}
\leq (\kappa-1)C^{\frac{1}{1-\kappa}}\kappa^{\frac{\kappa}{1-\kappa}}.$$

In order to get the lower bound, we introduce a small $\delta>0$. For the sake of clarity we call $\varepsilon:=\left(\frac{v}{u}\right)^{1/\kappa}.$ Note that for $\lambda>0$
\begin{multline*}E_{W}\left[e^{-\lambda^* H(v)}\right]=E_{W}\left[e^{-\lambda^* H(v)}{\bf 1}_{H(v)<({1-\delta})\varepsilon}\right]\\+E_{W}\left[e^{-\lambda^* H(v)}{\bf 1}_{(1-\delta)\varepsilon \leq H(v) \leq ({1+\delta})\varepsilon}\right]+E_{W}\left[e^{-\lambda^* H(v)}{\bf 1}_{H(v)>({1+\delta})\varepsilon}\right]\\:=J_1+J_2+J_3.\\
\end{multline*}
We are going to show that $J_1+J_3 \ll E_{W}\left[e^{-\lambda^* H(v)}\right]$.
We call $F(x)=P_{W}[H(v)<x]$.
By the Cramer-Chernoff inequality, for $x<\varepsilon,$  one gets 
\begin{multline}\label{cramer}F(x)\leq\exp\left(\lambda(x) x -C(1-o(1))v\lambda(x)^{\kappa}\right)\\=\exp\left(-C(1-o(1))v(1-\kappa)\lambda(x)^\kappa\right)\\
=\exp\left[-(1-o(1))C^{\frac{1}{1-\kappa}}(1-\kappa)\kappa^{\frac{\kappa}{1-\kappa}}v^{\frac{1}{1-\kappa}}x^{\frac{\kappa}{\kappa-1}}\right].
\end{multline}
Recall that 
$$E\left[e^{-\lambda^*H(v)}\right]=e^{-C(1+o(1))(C\kappa)^{\frac{\kappa}{1-\kappa}} u^{\frac{1}{1-\kappa}}}.$$
We deduce that for $\alpha=2(1-\kappa)^{\frac{\kappa-1}{\kappa}},$
$$F(\alpha\varepsilon)\ll E_{W}\left[e^{-\lambda^* H(v)}\right].$$ For this $\alpha$, we have

\begin{multline}\label{jun}
 J_1\leq F(\alpha\varepsilon)+\int_{\alpha\varepsilon}^{(1-\delta)\varepsilon} e^{-\lambda^* x} d F(x)\\=e^{-(1-\delta)\varepsilon \lambda^*} F((1-\delta)\varepsilon)+(1-e^{-\alpha\varepsilon})F(\alpha\varepsilon)+\lambda^*\int_{\alpha\varepsilon}^{(1-\delta)\varepsilon}e^{-\lambda^* x}F(x)d x.
\end{multline}

Our goal is to use (\ref{cramer}) in order to bound $F$ in the last equation. The problem is that the $o(1)$ in (\ref{cramer}) depends on $x$. We are going to use the monotonicity of $F(x)$ in order to get an uniform bound. Let $\eta<\delta/1000$, $n>\frac{\kappa}{\alpha(1-\kappa)\eta}$. For $1\leq k \leq n$, we set $x_k=k\varepsilon/n.$ 
Using (\ref{cramer}), there exists $v_0$ such that, for all $v>v_0$, and $1\leq k \leq n$, 
$$F(x_k)\leq \exp\left[-(1-\eta)C^{\frac{1}{1-\kappa}}(1-\kappa)\kappa^{\frac{\kappa}{1-\kappa}}v\left(\frac{x_k}{v}\right)^{\frac{\kappa}{\kappa-1}}\right].$$
  Note that for $x_{k-1}<x<x_k$, $x>\alpha\varepsilon$,
   and $v>v_0$,
$$F(x)\leq F(x_\kappa)\leq \exp\left[-(1-\eta)C^{\frac{1}{1-\kappa}}(1-\kappa)\kappa^{\frac{\kappa}{1-\kappa}}v^{\frac{1}{1-\kappa}}\left(x+\frac{\varepsilon}{n}\right)^{\frac{\kappa}{\kappa-1}}\right].$$
By the concavity of the function $x\rightarrow x^{\frac{\kappa}{\kappa-1}},$ and the condition $\varepsilon>x>\alpha\varepsilon$, we get easily
$$\left(x+\frac{\varepsilon}{n}\right)^{\frac{\kappa}{\kappa-1}}\geq x^{\frac{\kappa}{\kappa-1}}+\frac{1}{n}\frac{\kappa}{\alpha(\kappa-1)}\alpha^{\frac{\kappa}{\kappa-1}}\varepsilon^{\frac{\kappa}{\kappa-1}}\geq(1-\eta)x^{\frac{\kappa}{\kappa-1}}.$$
We deduce that for every $\varepsilon>x>\alpha\varepsilon$,
\begin{equation}\label{fonctionG}F(x)\leq \exp\left[-(1-\eta)^2 C^{\frac{1}{1-\kappa}}(1-\kappa)\kappa^{\frac{\kappa}{1-\kappa}}v^{\frac{1}{1-\kappa}}x^{\frac{\kappa}{\kappa-1}}\right]:=e^{G(x)}\end{equation}

Therefore, replacing $F$ by $e^{G}$ in (\ref{jun}), and doing the integration by parts in the other direction, we get
\begin{multline}J_1\leq e^{-(1-\delta)\varepsilon \lambda^*} e^{G((1-\delta)\varepsilon)}+(1-e^{-\alpha\varepsilon})e^{G(\alpha\varepsilon)}+\lambda^*\int_{\alpha\varepsilon}^{(1-\delta)\varepsilon}e^{-\lambda^* x}e^{G(x)}d x\\= 
e^{G(\alpha\varepsilon)}+\int_{\alpha\varepsilon}^{(1-\delta)\varepsilon}e^{-\lambda^* x}de^{G(x)}
.\end{multline}
Recalling the definition of $\alpha,$ $$ e^{G(\alpha\varepsilon)}\ll E\left[e^{-\lambda^*H(v)}\right], $$ and the integral can be bounded by
\begin{equation*}C' v^{\frac{\kappa}{1-\kappa}}\int_{\alpha\varepsilon}^{(1-\delta)\varepsilon}x^{\frac{1}{\kappa-1}}e^{-\lambda^* x}e^{G(x)}
d x.
\end{equation*}
Therefore, recalling (\ref{hdev}), and (\ref{lambda}) for estimates on $E_W\left[e^{-\lambda^*H(v)}\right]$, and the expressions of $\lambda(x)$ and $G$ respectively in (\ref{lambd}) and (\ref{fonctionG}), one gets
\begin{multline*}
 J_1 \left(E_W\left[e^{-\lambda^*H(v)}\right]\right)^{-1}\leq o(1)+P \sup_{x\in [\alpha \epsilon,(1-\delta) \epsilon]}
\exp\left(-C^{\frac{1}{1-\kappa}}\kappa^{\frac{\kappa}{1-\kappa}}v^{\frac{1}{1-\kappa}}\right.\\ \left.\left[(1-\eta)^2(1-\kappa)x^{\frac{\kappa}{1-\kappa}} +\kappa \varepsilon^{\frac{1}{\kappa-1}}-\varepsilon^{\frac{\kappa}{\kappa-1}}(1+o(1))\right]\right).\\
\end{multline*}
where P is some polynom in $(u,v)$ and the terms between the brackets come respectively from $e^G$, $e^{-\lambda x}$ and $\left(E_W\left[e^{-\lambda^*H(v)}\right]\right)^{-1}$. By a change of variable in the $\sup$, we get 
\begin{multline}\label{borneJ1}\frac{J_1}{E_W\left[e^{-\lambda^*H(v)}\right]}<o(1)+\\P\exp\left(-(Cv)^{\frac{1}{1-\kappa}}(\varepsilon\kappa)^{\frac{\kappa}{1-\kappa}}\inf_{s\in [\alpha,(1-\delta)]}\left[(1-\eta)^2s^{\frac{\kappa}{\kappa-1}}(1-\kappa)+\kappa s -1+o(1)\right]\right).\end{multline}

For $\eta$ and $o(1)$ very small, $$\inf_{s\in [\alpha,(1-\delta)]}\left[(1-\eta)^2s^{\frac{\kappa}{\kappa-1}}(1-\kappa)+\kappa s -1+o(1)\right]$$ is positive by concavity of the function $s\rightarrow s^{\frac{\kappa}{\kappa-1}}(1-\kappa)$ , therefore as an easy consequence

$$J_1\ll E_W\left[e^{-\lambda^*H(v)}\right].$$

We now deal with $J_3$.
As before we get
$$J_3<e^{-(1+\delta)\varepsilon \lambda^*} F((1+\delta)\varepsilon)+\lambda^*\int_{(1+\delta)\varepsilon}^{\infty}e^{-\lambda^* x}F(x)d x
$$
for $\beta >0$, as $F(x)\leq 1$
$$\lambda^*\int_{\beta\varepsilon}^{\infty}e^{-\lambda^* x}F(x)d x\leq e^{-\beta \lambda^*\varepsilon}=\exp{\left(-\beta(C\kappa)^{\frac{1}{1-\kappa}}u^{\frac{1}{1-\kappa}}\right)}$$
therefore for some  $\beta$ depending on $\kappa$,  
$$R(\varepsilon):=\lambda^*\int_{\beta\varepsilon}^{\infty}e^{-\lambda^* x}F(x)d x\ll E_{W}\left[e^{-\lambda^* H(v)}\right].$$
by the same argument as for $J_1$, we get that, for any $\varepsilon<x<\beta \varepsilon$, for $v$ large enough, 
$$F(x)\leq e^{G(x)};$$
therefore  
\begin{equation*} J_3-R(\varepsilon)\leq  e^{-(1+\delta)\varepsilon \lambda^*} e^{G((1+\delta)\varepsilon)}+\lambda^*\int_{(1+\delta)\varepsilon}^{\beta \varepsilon}e^{-\lambda^* x}e^G(x)d x 
\end{equation*}
By the same computation as we did to get to (\ref{borneJ1}), we have
\begin{multline}\label{borne J3}\frac{J_3}{E_W\left[e^{-\lambda^*H(v)}\right]}<o(1)+\\P\exp\left(-(Cv)^{\frac{1}{1-\kappa}}(\varepsilon\kappa)^{\frac{\kappa}{1-\kappa}}\inf_{s\in [(1+\delta,\beta]}\left[(1-\eta)^2s^{\frac{\kappa}{\kappa-1}}(1-\kappa)+\kappa s -1+o(1)\right]\right).\end{multline}
As before, we can take $\eta$ small and get $$J_3\ll  E_W\left[e^{-\lambda^*H(v)}\right].$$
Therefore we get that, as $v\rightarrow\infty,$
$$J_2>\frac{1}{2} E_W\left[e^{-\lambda^*H(v)}\right].$$
Recall that 
$$J_2=E_{W}\left[e^{-\lambda^* H(v)}{\bf 1}_{(1-\delta)\varepsilon \leq H(v) \leq ({1+\delta})\varepsilon}\right]\leq e^{-\lambda*(1-\delta)\varepsilon}P_W\left[H(v)<(1+\delta)\varepsilon\right].$$
Note that the preceding computations remain true for $u':=(1+\delta)^\kappa u,$ whence
$$\liminf_{v\rightarrow \infty} \frac{\log P_W\left[H(v)<\left(\frac{v}{u}\right)^{1/\kappa}\right]}{u^{\frac{1}{1-\kappa}}}>(1+\delta)^{\frac{\kappa}{1-\kappa}}
 ((1-\delta)\kappa-1)C^{\frac{1}{1-\kappa}}\kappa^{\frac{\kappa}{1-\kappa}}.$$
Taking the limit as $\delta \rightarrow 0,$ we get the result.
 
It remains to prove lemma \ref{equ}, which is the purpose of the next section.

\subsection{Proof of Lemma \ref{equ}.} 
Let $\nu =1 -2\kappa$, and
\begin{multline*}
D_\nu=\int_{0}^{\tau_r}
g^{-1}(\gamma(s))^{\nu}\exp{\left(-\frac{4}{g^{-1}(\gamma(s))}-8\lambda
  g^{-1}(\gamma(s))\right)}d s.\\=
\int_{-\infty}^{\infty}
g^{-1}(s)^{\nu}\exp{\left(-\frac{4}{g^{-1}(s)}-8\lambda
  g^{-1}(s)\right)}L_{\tau_r}^s d s.\\
=\left(\int_{-\infty}^0 + \int_{0}^{g(a)} + \int_{g(a)}^{\infty}\right)g^{-1}(s)^{\nu}\exp{\left(-\frac{4}{g^{-1}(s)}-8\lambda
  g^{-1}(s)\right)}L_{\tau_r}^s d s\\
:=I_1+I_2+I_3,
\end{multline*}
where $a$ is such that $a>1/\lambda$ and $$\frac{e^{4\lambda
    a}}{4\lambda a} = \log{\frac{1}{\lambda}}\log\log\log{\frac{1}{\lambda}}.$$

We shall use the following consequence of the law of large numbers :
let $f:\mathbb{R}\rightarrow\mathbb{R}$ such that
$\int_{\mathbb{R}}|f(x)|d x <\infty$, then
\begin{equation}\label{lln}
\lim_{r\rightarrow\infty}
\frac{1}{r}\int_{\mathbb{R}}f(x)L_{\tau_r}^x dx = \int_{\mathbb{R}}f(x)dx.
\end{equation}
Note that, for $x<1$ and $\lambda<1/4,$  $$|g(x)|=\int_x^1 \frac{e^{2/s
    + 4\lambda s}}{s^{1-\kappa}}d s\leq \frac{e^{2/x
    + 1}}{x^{1-\kappa}}, $$ therefore, for some constant $c>0,$ for
    all 
$x\leq0$ and $\lambda<1/4$ we have \begin{equation}\label{logx}\frac{2}{g^{-1}(x)}\geq \log\frac{|x|}{c}.\end{equation}
On the other hand
\begin{multline*}I_1=\int_{-\infty}^0 g^{-1}(s)^{\nu}\exp{\left(-\frac{4}{g^{-1}(s)}-8\lambda
  g^{-1}(s)\right)}L_{\tau_r}^s d s\\
\leq  \int_{-\infty}^0 g^{-1}(s)^{\nu}\exp{\left(-\frac{4}{g^{-1}(s)}\right)}L_{\tau_r}^s d s.
\end{multline*}
Using (\ref{logx}), it is not difficult to check that $g^{-1}(s)^{\nu}\exp{\left(-\frac{4}{g^{-1}(s)}\right)}$ is integrable on $(-\infty,0)$, therefore an application of (\ref{lln}) lays
$$I_1=O(r).$$
Let us now treat $I_3.$
Note that for $y\geq a$, $y\lambda\rightarrow \infty$ and for some constant $c>0$
$$\frac{1}{c}\int_1^y \frac{e^{4\lambda s}}{s^{1-\kappa}}d s
 \leq g(y)
\leq c \int_1^y \frac{e^{4\lambda s}}{s^{1-\kappa}}d s$$ and
\begin{equation}\label{gy}
 \int_1^y \frac{e^{4\lambda s}}{s^{1-\kappa}}d s=\frac{1}{(4\lambda)^\kappa} \int_{4\lambda}^{4\lambda y}
    \frac{e^{s}}{s^{1-\kappa}}d s=(1+o(1)) \left(\frac{e^{4\lambda y}}{4\lambda y}\right).
\end{equation}

As $y\lambda\rightarrow \infty$, we get $$g(y)\leq 2c e^{4\lambda y}.$$ Therefore, 
 for $x\geq g(a),$ 
  $2c e^{4\lambda g^{-1}(x)}\geq x$, so $g^{-1}(x)\geq\frac{1}{4\lambda}\log \frac{x}{2c}.$
Therefore, using (\ref{lln}) we get, for some constant $c'>0$,
\begin{multline*}I_3 \leq\int_{g(a)}^{\infty} (g^{-1}(x))^{\nu\vee 0}
  e^{-8\lambda g^{-1}(x)} L_{\tau_r}^x d x
\leq c'\int_{g(a)}^\infty \left(\frac{\log (x/2c)}{4\lambda}\right)^{\nu\vee 0}
x^{-2}  L_{\tau_r}^x d x\\
\leq c'(g(a))^{-1/2} \int_1^{\infty} \left(\frac{\log (x/2c)}{4\lambda}\right)^{\nu\vee 0}
x^{-3/2}  L_{\tau_r}^x d x= o\left(\frac{r}{\lambda^{\nu\vee 0}}\right).
\end{multline*}

To deal with $I_2$, note that, by the definition of $a$ and (\ref{gy}),
$$r\gg g(a) \log \log g(a).$$
Therefore we can apply statement \ref{tl} to get
$$I_2=r(1+o(1))\int_0^{g(a)}g^{-1}(s)^{\nu}\exp{\left(-\frac{4}{g^{-1}(s)}-8\lambda
  g^{-1}(s)\right)}d s.$$
By a change of variables $g^{-1}(s)=y$, as $\lambda\rightarrow 0$,
  the last integral is equal to
$$
\int_1^a \frac{e^{-\frac{2}{y}-4\lambda
    y}}{y^{1-\kappa-\nu}}d y=(1+o(1))\frac{1}{(4\lambda)^{\kappa+\nu}}\int_{4\lambda}^{4\lambda
    a} \frac{e^{-u}}{u^{1-(\nu+\kappa)}}d u.
$$
Recalling the definition of $\nu$ we have $\nu+\kappa
=1-\kappa>0$, then
$$I_2=r(1+o(1))\frac{\Gamma(1- \kappa)}{(4\lambda)^{1-\kappa}}.$$

This finishes the proof of the first part of lemma \ref{equ}, as $1-\kappa>\nu\vee 0$.
\\
To treat the case $\nu=-2\kappa$, let $b<1$ be such that $b\rightarrow 0$ and
$-g(b)=o\left(\frac{r}{\log\log r}\right).$
As before, we separate the integral as follows
\begin{multline*}
D_\nu
=\left(\int_{-\infty}^{g(b)} + \int_{g(b)}^{g(a)} + \int_{g(a)}^{\infty}\right)g^{-1}(s)^{\nu}\exp{\left(-\frac{4}{g^{-1}(s)}-8\lambda
  g^{-1}(s)\right)}L_{\tau_r}^s d s\\
:=I_1'+I_2'+I_3'.
\end{multline*}
$I_3'$ is similar to the precedent case, with $\nu<0$, so we get $I_3'=o(r)$. 
We have easily 
$$I'_1\leq e^{-\frac{1}{b}}\int_{-\infty}^0 g^{-1}(s)^{\nu}\exp{\left(-\frac{3}{g^{-1}(s)}-8\lambda
  g^{-1}(s)\right)}L_{\tau_r}^s d s.$$ The integral is a $O(r)$ by the same proof as for $I_1$, therefore$$I'_1=o(r).$$

By the same proof as for $I_2$, we
get
$$I'_2=r(1+o(1))I''_2,$$ with
$$I''_2=\int_b^a \frac{e^{-\frac{2}{y}-4\lambda
    y}}{y^{1+\kappa}}d y=\int_b^1 \frac{e^{-\frac{2}{y}-4\lambda
    y}}{y^{1+\kappa}}d y+\int_1^a \frac{e^{-\frac{2}{y}-4\lambda
    y}}{y^{1+\kappa}}d y.$$
The first part converges, by dominated convergence, to 
$$D:=\int_0^1 \frac{e^{-\frac{2}{y}}}{y^{1+\kappa}}d y,$$
 and the second part is equal to
$$ (4\lambda)^{\kappa}\int_{4\lambda}^{4\lambda
    a} \frac{e^{-8\lambda/u-u}}{u^{1+\kappa}}d u.$$
One can easily check that the integral is bounded, therefore this part goes to zero. This finishes the proof of lemma \ref{equ}.
\subsection{Quenched Speedup for the diffusion.}
In this section we prove Theorem \ref{qs}.
The upper bound is a trivial consequence of Theorem \ref{pqs} , since
$$P_{W}[X_t>t^\kappa u]\leq P_W[H(t^\kappa u)<t].$$
To get the lower bound, let $\varepsilon>0.$ Note that
$$P_{W}[X_t>t^\kappa u]\geq P_{W}[H((1+\varepsilon)t^\kappa u)<t]P_W^{(1+\varepsilon)t^\kappa u}[H(t^{\kappa} u)>t].$$

Note that almost surely, for $t$ large enough, we can find 
$t^\kappa u <b<c< (1+\varepsilon)t^\kappa u$ such that
$$W_{\kappa}(b)-W_{\kappa}(c)> \frac{\varepsilon \kappa}{2}t^\kappa u.$$
It is clear that 
$$P_W^{2t^\kappa u}[H(t^{\kappa} u)>t]\geq P_W^c[H(b)>t].$$
By the same computations as in \ref{lb}, one gets easily that $$P_W^c[H(b)>t]>1/2$$ for $t$ large enough. Taking the limit as $\varepsilon\rightarrow 0$, this finishes the proof of Theorem \ref{qs}.

{\bf Acknowledgement :} We are thankful to Yueyun Hu for pointing this subject to us and for many helpful discussions, and to an anonymous referee for detailed and constructive remarks on the manuscript.

\bibliographystyle{plain}
\bibliography{biblio}

\begin{thebibliography}{10}

\bibitem{Bertoin:1996}
J.~Bertoin.
\newblock {\em L\'evy processes}.
\newblock Cambridge Tracts in Mathematics, 1996.

\bibitem{Bobkov:2002}
S.G. Bobkov and F.~G\"otze.
\newblock Muckenhoupt's condition via \text{R}iccati and
  \text{S}turm-\text{L}iouville equations.
\newblock {\em Pr\'eprint}, 2002.

\bibitem{Borodin:2002}
A.N. Borodin and P.~Salminen.
\newblock {\em Handbook of Brownian Motion}.
\newblock Probability and its application. Birkha\"user, 2002.

\bibitem{brox:1986}
T.~Brox.
\newblock A one-dimensionnal diffusion process in a \text{W}iener medium.
\newblock {\em The Annals of Probability}, 14:1206--1218, 1986.

\bibitem{chernov1967rmc}
A.A. Chernov.
\newblock {Replication of a multicomponent chain, by the lightning mechanism}.
\newblock {\em Biophysics}, 12(2):336--341, 1967.

\bibitem{Csorgo:1981}
M.~Cs\"org\"o and P.~R\'ev\'esz.
\newblock {\em Strong approximations in Probability and Statistics}.
\newblock Academic Press, New York, 1981.

\bibitem{Gantert:2008}
A~Fribergh, N.~Gantert, and S.~Popov.
\newblock On the slowdown and speedup of transient random walks in random
  environment.
\newblock {\em Preprint}, 2008.

\bibitem{Hu:2004}
Y.~Hu and Z.~Shi.
\newblock Moderate deviations for diffusions with brownian potentials.
\newblock {\em The Annals of Probability}, 32(4):3191--3220, 2004.

\bibitem{Hu:1999}
Y.~Hu, Z.~Shi, and M.~Yor.
\newblock Rates of convergence of diffusions with drifted brownian potentials.
\newblock {\em Trans. Amer. Math. Soc.}, 351(26):3915--3934, 1999.

\bibitem{ito1974dpa}
K.~It{\=o} and H.P. McKean.
\newblock {\em {Diffusion processes and their sample paths}}.
\newblock Springer, 1974.

\bibitem{Jeulin:1981}
T.~Jeulin and M.~Yor.
\newblock Sur les distributions de certaines fonctionnelles du mouvement
  brownien.
\newblock In {\em S\'em. Prob. XV}, number 850 in Lect. Notes in Math., pages
  210--226. Springer, 1981.

\bibitem{Kawazu:1993}
K.~Kawazu and H.~Tanaka.
\newblock On the maximun of a diffusion process in a drifted brownian
  environment.
\newblock In {\em S\'eminaires de probabilit\'es XXVII}, number 1557 in Lecture
  Notes in Mathematics, pages 78--85. Springer, 1993.

\bibitem{Kawazu:1997}
K.~Kawazu and H.~Tanaka.
\newblock A diffusion process in a brownian random environment with drift.
\newblock {\em J. Math. Soc. Japan}, 49:189--211, 1997.

\bibitem{kesten1975llr}
H.~Kesten, M.V. Kozlov, and F.~Spitzer.
\newblock {A limit law for random walk in a random environment}.
\newblock {\em Compositio Math}, 30:145--168, 1975.

\bibitem{Pitman:1982}
J.~Pitman and M.~Yor.
\newblock {Sur une decomposition des ponts de Bessel}.
\newblock {\em Functional Analysis in Markov Processes}, pages 276--285, 1982.

\bibitem{Revuz:1994}
D.~Revuz and M.~Yor.
\newblock {\em Continuous Martingales ans Brownian Motion}.
\newblock Springer, Berlin, 1994.

\bibitem{Rogers:1986}
L.C.G. Rogers and D.~Williams.
\newblock {\em Diffusions, Markov processes and Martingales}, volume~2.
\newblock Wiley series in probability and mathematical statistics, 1986.

\bibitem{Salminen:2001}
P.~Salminen and I.~Norros.
\newblock On busy periods of the unbounded brownian storage.
\newblock {\em Queuing Systems}, 39:317--333, 2001.

\bibitem{schumacher1985drc}
S.~Schumacher.
\newblock {Diffusions with random coefficients}.
\newblock {\em Particle systems, random media and large deviations}, pages
  351--356, 1985.

\bibitem{sinai1983lbo}
Y.G. Sinai.
\newblock {The limiting behavior of a one-dimensional random walk in a random
  medium}.
\newblock {\em Theory of Probability and its Applications}, 27:256, 1983.

\bibitem{solomon1975rwr}
F.~Solomon.
\newblock {Random walks in a random environment}.
\newblock {\em The annals of probability}, pages 1--31, 1975.

\bibitem{Talet:2007}
M.~Talet.
\newblock Annealed tail estimates for a brownian motion in a drifted brownian
  potential.
\newblock {\em The Annals of Probability}, 35(1):32--67, 2007.

\bibitem{temkin1972odr}
D.E. Temkin.
\newblock {One-dimensional random walks in a two-component chain}.
\newblock In {\em Soviet Math. Dokl}, volume~13, pages 1172--1176, 1972.

\bibitem{watanabe1977comparison}
N.I.A. Watanabe.
\newblock {A comparison theorem for solutions of stochastic differential
  equations and its applications}.
\newblock {\em Osaka J. Math}, 14:619--633, 1977.

\bibitem{Werner:1995}
W.~Werner.
\newblock Some remarks on perturbed reflecting brownian motion.
\newblock {\em Seminaire de Probabilit\'es (Strasbourg)}, 29:37--43, 1995.

\end{thebibliography}
\end{document}